\documentclass[12pt]{article}

\usepackage{amsmath}
\usepackage{amssymb}
\usepackage{array}
\usepackage{delarray}
\usepackage{amsfonts}
\usepackage{amsthm}
\usepackage[english]{babel}
\usepackage{enumerate}
\usepackage{fancyhdr}
\usepackage{slashbox}
\usepackage{graphicx}
\usepackage{float}
\usepackage{epstopdf}
\usepackage{amsmath,amsthm,amsopn,amstext,amscd,amsfonts,amssymb,amsbsy}
\usepackage{mathrsfs}

\RequirePackage[OT1]{fontenc}

\theoremstyle{plain}
\newtheorem{theorem}{Theorem}
\newtheorem{coro}{Corollary}
\newtheorem{lem}{Lemma}
\newtheorem{pro}{Proposition}

\newtheorem{definition}{Definition}

\newcommand{\Bea} {\begin{eqnarray*}}
\newcommand{\Eea} {\end{eqnarray*}}

\newcommand{\E}{\mathbb{E}}
\newcommand{\p}{\mathbb{P}}

\newcommand{\eps}{\varepsilon}

\begin{document}
\title{Regression in random design and Bayesian warped wavelets estimators}
\author{Thanh Mai Pham Ngoc \\ \it Universit{\'e} Pierre et Marie Curie, Paris VI}
\date{}
\maketitle
\begin{abstract}
In this paper we deal with the regression problem in a random design
setting. We investigate asymptotic optimality under minimax point of
view of various Bayesian rules based on warped wavelets and show
that they nearly attain optimal minimax rates of convergence over
the Besov smoothness class considered. Warped wavelets have been
introduced recently, they offer very good computable and
easy-to-implement properties while being well adapted to the
statistical problem at hand. We particularly put emphasis on
 Bayesian rules leaning on small and large variance Gaussian priors and discuss their simulation performances comparing them with a hard
thresholding procedure.
\end{abstract}
 \vspace{9pt} \noindent {\bf Key words and phrases:}
{nonparametric regression}, {random design}, {warped wavelets},
{Bayesian methods}
\par \vspace{9pt} \noindent {\bf MSC 2000 Subject Classification}
62G05 62G08 62G20 62C10
\par

\section{Introduction}

We observe independant pairs of variables $(X_{i},Y_{i})$, for
$i=1,\dots,n$, under a random design regression model:
\begin{equation}\label{modele_de_base}
Y_{i}=f(X_{i})+\eps_{i}, \quad 1\leq i \leq n,
\end{equation}
where $f$ is an unknown regression function that we aim at
estimating, and $\eps_{i}$ are independant normal errors with
$\E(\eps_{i})=0$, $\textrm{Var}(\eps_{i})=\sigma^2<\infty$. The
design points $X_{i}$ are assumed to be supported in the interval
$[0,1]$ and have a density $g$ which will be supposed to be known.
Furthermore we assume that the design density $g$ is bounded from
below, i.e. $0<m\leq g$, where $m$ is a constant. Many approaches
have been proposed to tackle the problem of regression in random
design, we mention among others the work of Hall and Turlach
\cite{HallTurlach},
 Kovac and Silverman \cite{ KovacSilverman}, Antoniadis et \textit{al}. \cite{
 AntoniaGregoireVial}, Cai and Brown \cite{ Cai_Brown_annals} and the model selection point
of view adopted by Baraud \cite{Baraud}. \\The present paper
provides a Bayesian approach to this problem based on
\textit{warped} wavelet basis. Warped wavelets basis
$\{\psi_{jk}(G)\; j\geq -1,k\in \mathbb{Z}\}$ in regression with
random design were recently introduced by Kerkyacharian and Picard
in \cite{Kerk2}. The authors proposed an approach which would depart
as little as possible from standard wavelet thresholding procedures
which enjoy optimality and adaptivity properties. These procedures
have been largely investigated in the case of equispaced samples
(see a series of pioneered articles by Donoho \textit{et al.}
\cite{DonohoAdaptingUnknown}, \cite{DonohoAsymptotia},
\cite{DonohoIdealSpatial}). Kerkyacharian and Picard actually
pointed out that expanding the unknown regression function $f$ in
the warped basis instead of the standard wavelets basis could be
very interesting. Of course, this basis has no longer the
orthonormality property nonetheless it behaves under some conditions
as standard wavelets. Kerkyacharian and Picard investigated the
properties of this new basis and showed that not only is it well
adapted to the statistical problem at hand by avoiding unnecessary
calculations but it also offers very good theoretical features while
being easily implemented. More recently Brutti \cite{Brutti}
highlighted their easy-to-implement computational properties. \\
The novelty of our contribution lies in the combination of Bayesian
techniques and \textit{warped} wavelets to treat regression in
random design. We actually want to investigate whether this yields
optimal theoretical results and promising pratical performances,
which will prove to be the case. We do not deal with the case of an
unknown design density $g$ which requires further machinery and will
be the object of another paper.\\Bayesian techniques for shrinking
wavelet coefficients have become very popular in the last few years.
The majority of them were devoted to fixed design regression scheme.
Let us cite among others, papers of Abramovich \textit{et al.}
\cite{Abra1}, \cite{Abra2}, Clyde et al. \cite{ ClydeFlexible},
\cite{ClydeMultiple}, \cite{ClydeEmpirical}, \cite{Autin}, Chipman
\textit{et al.} \cite{Chipman}, Rivoirard
\cite{RivoirardBayesianMMS}, Pensky \cite{PenskyFrequentist} in the
case of i.i.d errors not necessarily Gaussian.  \\
Most of those works are taking as distribution prior a mixture of
Gaussian distributions. In particular, Abramovich \textit{et al.} in
\cite{Abra1} and \cite{Abra2} have explored optimality properties of
Gaussian prior mixed with a point mass at zero and which may be
viewed as an extreme case of a Gaussian mixture:
$$ \beta_{jk} \sim \pi_{j}N(0,\tau_{j}^{2})+ (1-\pi_{j})\delta(0)$$
where $\beta_{jk}$ are the wavelet coefficients of the unknown
regression function, $\tau_{j}^{2}=c_{1}2^{-j\alpha}$ and
$\pi_{j}=\min(1,c_{2}2^{-j\beta})$ are the hyperparameters. This
particular form was devised to capture the sparsity of the expansion
of the signal in the wavelets basis.
\\
Our approach will consist in a first time in using the same prior
but in the context of warped wavelets. In Theorem 1 we show that the
Bayesian estimator built using warped wavelets with this prior and
this form of hyperparameters achieves the optimal minimax rate
within logarithmic term on the considered Besov functional space.
Unfortunately, the Bayesian estimator turns out not to be adaptive.
Indeed, the hyperparameters depend on the Besov smoothness class
index. In order to compensate this drawback, Autin \textit{et al.}
in \cite{Autin} suggested to consider Bayesian procedures based on
Gaussian prior with large variance. Following this suggestion, we
will consider priors still specified in terms of a normal density
mixed with a point mass at zero but with large variance Gaussian
densities. In Theorem 2 we prove again that the Bayesian estimator
built with this latter form of prior, still combined with warped
wavelets achieves nearly optimal minimax rate of convergence while
being adaptive. Eventually, our simulations results highlight the
very good performances and behaviour of these Bayesian procedures
whatever the regularity of the test functions, the noise level and
the design
density which can be far from the uniform case may be.\\
This paper is organized as follows. In section 2 some necessary
methodology is given: we start with a short review of wavelets and
warped wavelets, explain the prior model and discuss the two
hyperparameters form we consider. We give in section 3 some
definitions of functional spaces we consider. In section 4, we
investigate the performances of our Bayesian estimators in terms of
minimax rates in two cases: the first one when the Gaussian prior
has small variance, the second case focuses on Gaussian prior with
large variance. Section 5 is devoted to simulation results and
discussion. Finally, all proofs of main results are given in the
Appendix.

\section{Methodology}
\subsection{Warped bases}
Wavelet series are generated by dilations and translations of a
function $\psi$ called the mother wavelet. Let $\phi$ denote the
orthogonal father wavelet function. The function $\phi$ and $\psi$
are compactly supported. Assume $\psi$ has $r$ vanishing moments.
Let:
$$\phi_{jk}(x)=2^{j/2}\phi(2^jx-k) , \quad j,k \in \mathbb{Z}$$
$$\psi_{j,k}(x)=2^{j/2}\psi(2^jx-k), \quad j,k \in \mathbb{Z}.$$
For a given square-integrable function $f$ in $\mathbb{L}_2[0,1]$,
let us denote
$$\zeta_{j,k}=<f,\psi_{j,k}>.$$
In this paper, we use decompositions of 1- periodic functions on
wavelet basis of $\mathbb{L}_2[0,1]$. We consider periodic
orthonormal wavelet bases on $[0,1]$ which allow to have the
following series representation of a function $f$ :
\begin{equation} \label{decomposition_wave} f(x)=\sum_{j\geq-1}\sum_{k=0}^{2^j-1}\zeta_{jk}\psi_{jk}(x)\end{equation}
where we have denoted $\psi_{-1,k}=\phi_{0,k}$ the scaling function.
\\
We are now going to give the essential background of \textit{warped}
wavelets which were introduced in details in \cite{Kerk2}. First of
all let us define
\begin{equation}
G(x)=\int_{0}^{x}g(x)dx.
\end{equation}
G is assumed to be a known function, continuous and strictly
monotone from $[0,1]$ to $[0,1]$. \\
Let us expand the regression function $f$ in the following sense:
$$f(G^{-1})(x)=\sum_{j\geq-1}\sum_{k=0}^{2^j-1}\beta_{jk}\psi_{jk}(x)$$
or equivalently
$$f(x)=\sum_{j\geq-1}\sum_{k=0}^{2^j-1}\beta_{jk}\psi_{jk}(G(x))$$
where
$$ \beta_{jk}=\int f(G^{-1})(x)\psi_{jk}(x)dx=\int
f(x)\psi_{jk}(G(x))g(x)dx.$$ \\
Hence one immediately notices that expanding $f(G^{-1})$ in the
standard basis is equivalent to expand $f$ in the new
\textit{warped} wavelets basis $\{\psi_{jk}(G),j\geq -1, k\in
\mathbb{Z}\}$. This may give a natural explanation that in the
follow-on, regularity conditions will be expressed not for $f$ but
for $f(G^{-1})$.\\

We set
$\hat{\beta}_{jk}=(1/n)\sum_{i=1}^{n}\psi_{jk}(G(X_{i}))Y_{i}$.
$\hat{\beta}_{jk}$ is an unbiased estimate of $\beta_{jk}$ since
\begin{eqnarray*}
 \E(\hat{\beta}_{jk})&=&(1/n)\sum_{i=1}^{n}\E(\psi_{j,k}(G(X_i))(f(X_i)+\epsilon_i))=\E(\psi_{j,k}(G(X))f(X)\\
 &=&\int f(x)\psi_{jk}(G(x))g(x)dx= \int
 f(G^{-1})(x)\psi_{jk}(x)dx=\beta_{jk}.
 \end{eqnarray*}

 \subsection{Priors and estimators}
We set in the following
\begin{equation}\label{bruit_gamma}
\gamma_{jk}^{2}=
\frac{\sigma^2}{n^{2}}\sum_{i=1}^{n}\psi_{jk}^{2}(G(X_{i})) .
\end{equation}
\\ As in Abramovich \textit{et al.}
(see \cite{Abra1}, \cite{Abra2}), we use the following prior on the
wavelet coefficients $\beta_{jk}$ of the unknown function $f$ with
respect to the \textit{warped} basis $\{\psi_{jk}(G),j\geq -1, k\in
\mathbb{Z}\}$:
$$\beta_{jk}\sim \pi_{j}N(0,\tau_{j}^{2})+(1-\pi_{j})\delta(0).$$
Considering the $\mathbb{L}_{1}$ loss, from this form of prior we
derive the following Bayesian rule which is the posterior median:
\begin{equation}\label{bayesien_beta}
\tilde{\beta}_{jk}=Med(\beta_{jk}|\hat{\beta}_{jk})=\mbox{sign}(\hat{\beta}_{jk})\max(0,\zeta_{jk})
\end{equation}
where
\begin{equation}\zeta_{jk}=\frac{\tau_{j}^{2}}{\gamma_{jk}^{2}+\tau_{j}^{2}}|\hat{\beta}_{jk}|-
\frac{\tau_{j}\gamma_{jk}}{\sqrt{\gamma_{jk}^{2}+\tau_{j}^{2}}}\Phi^{-1}\bigg(\frac{1+\min(\eta_{jk},1)}{2}\bigg)
\end{equation}
where $\Phi$ is the normal cumulative distributive function and
\begin{equation}
\eta_{jk}=\frac{1-\pi_{j}}{\pi_{j}}\frac{\sqrt{\tau_{j}^{2}+\gamma_{jk}^{2}}}{\gamma_{jk}}
\exp\bigg(-\frac{\tau_{j}^{2}\hat{\beta}_{jk}^{2}}{2\gamma_{jk}^{2}(\tau_{j}^{2}+\gamma_{jk}^{2})}\bigg).
\end{equation}
We set :
\begin{equation}
w_{j}(n):=\frac{\pi_{j}}{1-\pi_{j}}.
\end{equation}
We introduce now the estimator of the unknown regression $f$
\begin{equation}\label{estimator}
\tilde{f}(x)=\sum_{j\leq
J}\sum_{k=0}^{2^j-1}\tilde{\beta}_{jk}\psi_{jk}(G(x))
\end{equation}
where $J$ is a parameter which will be precised later.\\
Note that in our case, the estimator resembles the usual ones in
\cite{Autin}, \cite{Abra1} and \cite{Abra2}, except that the
deterministic noise variance has been replaced by a stochastic noise
level $\gamma_{jk}^{2}$. Its expression is given by
(\ref{bruit_gamma}). This change will have a marked impact both on
the proofs of theorems by using now large deviations inequalities
and on simulations results.
\\
Futhermore, such $\mathbb{L}_{1}$ rule is of thresholding type.
Indeed, as underlined in \cite{Abra1} and \cite{Abra2},
$\tilde{\beta}_{jk}$ is null whenever $\hat{\beta}_{jk}$ falls below
a certain threshold $\lambda_{B}$. Some properties of the threshold
$\lambda_{B}$ that will be used in the sequel are given in
 lemma 1 in Appendix.

\subsubsection{Gaussian priors with small variance}
In this paper, two cases of hyperparameters will be considered. The
first one involves Gaussian priors with small variances. We will
state as suggested in Abramovich \textit{et al} (see \cite{Abra1},
\cite{Abra2}) :
\begin{equation} \label{hyperparametresmall}
\tau_{j}^{2}=c_{1}2^{-j\alpha} \quad
\pi_{j}=\min(1,c_{2}2^{-j\beta}),
\end{equation}
where $\alpha$ and $\beta$ are non-negative constants, $c_{1}, c_{2}
>0$.\\
This choice of hyperparameters is exhaustively discussed in
Abramovich \textit{et al.} \cite{Abra2}. The authors stressed that
this form of hyperparameters was actually designed in order to
capture
 the sparsity of wavelet expansion. They pointed out the connection between Besov
 spaces parameters and this particular form of hyperparameters. They
 investigate various practical choices.\\
 For this case of hyperparameters (\ref{hyperparametresmall}), the estimator of $f$ will be denoted $\hat{f}$.

\subsubsection{Gaussian priors with large variance}
 The second form of hyperparameters considered in the paper involves
Gaussian priors with large variance as suggested in Autin \textit{et
al}.
\cite{Autin}.\\
As a matter of fact, we suppose that the hyperparameters do not
depend on $j$ and we set :
\begin{equation}\label{hyperparametrebigtau}
\tau_{j}^{2}:=\tau(n)^{2}=1/\sqrt{n\log(n)}.
\end{equation}
Besides, $w_{j}(n):=w(n)$. We suppose that there exist $q_{1}$ and
$q_{2}$ such that for $n$ large enough
\begin{equation}\label{hyperparametrebigw}
{n}^{-q_{1}/2}\leq w(n)\leq {n}^{-q_{2}/2}.
\end{equation}
This form of hyperparameters was emphasized in \cite{Autin} in order
to mimic heavy tailed priors such as Laplace or Cauchy
distributions. Indeed, Johnstone and Silverman in
\cite{Johnstoneempiri1}, \cite{Johnstoneempiri2} showed that their
empirical Bayes approach for regular regression setting with a prior
mixing a heavy-tailed density and a point mass at zero proved
fruitful both in theory and practice. Pensky in
\cite{PenskyFrequentist} also underlined the efficiency of this
kind of hyperparameters.\\
 We underscore that contrary to the first form of
hyperparameters (\ref{hyperparametresmall}), this latter forms
(\ref{hyperparametrebigtau}) and (\ref{hyperparametrebigw}) lead to
an adaptive Bayesian estimator. \\For this case of hyperparameters
(\ref{hyperparametrebigtau}) and (\ref{hyperparametrebigw}), the
estimator of $f$ will be denoted $\check{f}$.

\vspace{6mm}

\section{Functional spaces}
In this paper, functional classes of interest are Besov bodies and
weak Besov bodies. Let us define them. Using the decomposition
(\ref{decomposition_wave}), we characterize Besov spaces by using
the following norm
$$\|f\|_{spq}=\left\{
\begin{array}{ll}
{\big[\sum_{j\geq-1}2^{jq(s+1/2-1/p)}\|(\beta_{j,k})_{k}\|^{q}_{\ell_{p}}\big]}^{1/q}
&
\mbox{if } q<\infty \\
\sup_{j\geq -1} 2^{j(s+1/2-1/p)}\|(\beta_{j,k})_{k}\|_{\ell_{p}} &
\mbox{if } q=\infty. \end{array} \right. $$ If $\max(0,1/p-1/2)<s<r$
and $p,q\geq 1$
$$ \quad f \in B^{s}_{p,q} \Longleftrightarrow \|f\|_{spq}<\infty.$$
The Besov spaces have the following simple relationship
$$ B^{s_{1}}_{p,q_{1}} \subset B^{s}_{p,q},\quad  \mbox{for } s_{1}>s \; \mbox{or } s_{1}=s \; \mbox{and  } q_{1}\leq q $$
and
$$ B^{s}_{p,q} \subset B^{s_{1}}_{p_{1},q}, \quad \mbox{for } p_{1}>p \; \mbox{and  } s_{1}\geq s-1/p+1/p_{1}
.$$ The index $s$ indicates the smoothness of the function. The
Besov spaces capture a variety of smoothness features in a function
including spatially inhomogeneous behavior when $p<2$.\\
We recall and stress that in this paper as mentioned above, the
regularity conditions will be expressed for the function $f(G^{-1})$
due to the \textit{warped} basis context. \\More precisely we shall
focus on the space $B^{s}_{2,\infty}$. We have in particular
\begin{equation} \label{besov} f\in B^{s}_{2,\infty}
\Longleftrightarrow \sup_{J\geq -1}2^{2Js}\sum_{j\geq J}\sum_{k}
\beta_{jk}^{2}< \infty.\end{equation} \\We define the Besov ball of
some radius R as $B^{s}_{2,\infty}(R)=\{f:\|f\|_{s2\infty}\leq R
\}$.
\\Let us define now the weak Besov space $W(r,2)$
\begin{definition}
Let $0<r<2$. We say that a function $f$ belongs to the weak Besov
body $W(r,2)$ if and only if:
\begin{equation}
\|f\|_{W_{r}}:=[\sup_{\lambda>0} \lambda^{r-2}
\sum_{j\geq-1}\sum_{k} \beta_{jk}^{2} I\{|\beta_{jk}\leq
\lambda|\}]^{1/2} < \infty.
\end{equation}
\end{definition}
And we have the following proposition
\begin{pro}
Let $0<r<2$ and $f\in W(r,2)$. Then
\begin{equation}\label{weakbesov}
\sup_{\lambda>0} \lambda^{r}\sum_{j\geq
-1}\sum_{k}I\{|\beta_{jk}|>\lambda\} \leq
\frac{2^{2-r}\|f\|_{W_{r}}^{2}}{1-2^{-r}}.
\end{equation}
\end{pro}
For the proof of this proposition see for instance \cite{Kerk3} . \\
To conclude this section, we have the following embedding
$$B^{s}_{2,\infty} \subset W_{2,2/(1+2s)}$$
which is not difficult to prove (see for instance \cite{Kerk3}).

\section{Minimax performances of the procedures}

\subsection{Bayesian estimators based on Gaussian priors with small
variances}

\begin{theorem}\label{Theo_1}
Assume that we observe model (\ref{modele_de_base}). We consider the
hyperparameters defined by (\ref{hyperparametresmall}). Set $J:=J_{\alpha}$ such that $2^{J_{\alpha}}=(3/(2n))^{-1/\alpha}$. \\
Let $\alpha>1$ and $\alpha\geq s$, then we have the following upper
bound:
\begin{equation}\label{borne_sup_1}
\sup_{f(G^{-1})\in
B^{s}_{2,\infty}(R)}\mathbb{E}\|\hat{f}-f\|^{2}_{2} =
\mathcal{O}((1/n)^{1-1/\alpha}\log^{2}(n))+\mathcal{O}((1/n)^{2s/\alpha}).
\end{equation}

\end{theorem}
\vspace{6mm} The optimal choice of the hyperparameter $\alpha$ in
Theorem \ref{Theo_1} should minimize the upper bound derived in
(\ref{borne_sup_1}). Consequently, let us choose now in
(\ref{borne_sup_1}) $\alpha=2s+1$, we immediately deduce the
following corollary.

\begin{coro}
If one chooses for the prior parameter $\alpha=2s+1$, one gets
$$ \sup_{f(G^{-1})\in B^{s}_{2,\infty}(R)}\E\|\hat{f}-f\|^{2}_{2}= \mathcal{O}(\log^{2}(n)n^{-2s/(2s+1)}).$$
\end{coro}
This corollary shows that with this specific choice of
hyperparameter $\alpha$, one recovers the minimax rate of
convergence up to a logarithmic factor that one achieves in a
uniform design.

\subsection{Bayesian estimators based on Gaussian priors with large
variance}
\begin{theorem}{\label{Theorem_2}}
We consider the model (\ref{modele_de_base}). We assume that the
hyperparameters are defined by (\ref{hyperparametrebigtau}) and
(\ref{hyperparametrebigw}). Set $J:=J_{n}$ such that
$2^{J_{n}}=n/\log n$, then we have :
$$\sup_{f(G^{-1})\in B^{s}_{2,\infty}(R)}\mathbb{E}\|\check{f}-f\|^{2}_{2}\leq C\left(\frac{\log(n)}{n}\right)^{2s/(2s+1)}.$$
\end{theorem}

It is worthwhile to make some comments about the results of Theorem
2. Here, the estimator turns out to be adaptive and contrary to the
similar results in Proposition 2 in \cite{Kerk2} we no longer have
the limitation on the regularity index $s>1/2$. Moreover,
Kerkyacharian and Picard \cite{Kerk2} had to stop
 the highest level $J$ such that $2^J=(n/\log(n))^{1/2}$, here we
 stop at the usual level $J_{n}$ such that $2^{J_{n}}=n/\log(n)$ one gets in standard thresholding .

\section{Simulations and discussion}
\vspace{6 mm} A simulation study is conducted in order to compare
the numerical performances of the two Bayesian estimators based on
warped wavelets and on Gaussian prior with small or large variance,
described respectively in section 2.2.1 and 2.2.2 and the hard
thresholding procedure using the universal threshold
$\sigma\sqrt{2\log(n)}$ based on warped basis introduced by
Kerkyacharian and Picard \cite{Kerk2} for the nonparametric
regression model in a random design setting. For more details on
Kerkyacharian and Picard procedure, the readers are referred to
Willer \cite{Willer_these}, see also \cite{Gannaz_these}. In fact,
we have decided to concentrate on the procedure of Kerkyacharian and
Picard because it is interesting to point out differences and
compare performances obtained by Bayesian procedures which apply
local thresholds and a universal threshold procedure.
\\The main difficulties lie in implementing the Bayesian procedures
with the stochastic variance (\ref{bruit_gamma}). Note also the
responses proposed by Amato et al. \cite{AmatoAntoniadisPensky} and
 Kovac and Silverman \cite{KovacSilverman}.
\\All the simulations done in the present paper have
been conducted with MATLAB and the Wavelet toolbox of MATLAB.\\
We consider here four test functions of Donoho and Johnstone
\cite{DonohoIdealSpatial} representing different level of spatial
variability. The test functions are plotted in Fig. 1. For each of
the four objects under study, we compare the three estimators at two
noise levels, one with signal-to-noise ratio $RSNR=4$ and another
with $RSNR=7$. As in Willer \cite{Willer_these} we also consider
different cases of design density which are plotted in Fig. 2. The
first two densities are uniform or slightly varying whereas
 the last two ones aim at depicting the case where a hole occurs in the density design.
 The sample size is equal to $n=1024$ and the wavelet we used is the Symmlet8.\\
In order to compare the behaviors of the estimators, the RMSE
criterion was retained. More precisely, if  $\hat{f}(X_{i})$ is the
estimated function value at $X_i$
 and $n$ is the sample size, then
 \begin{equation}\label{RMSE} RMSE=\sqrt{\frac{1}{n}\sum_{i=1}^{n}(\hat{f}(X_i)-f(X_{i}))^{2}}.\end{equation}\\
 The RMSE displayed in Tab. 1 are computed as the average over $100$ runs of expression (\ref{RMSE}). In each run,
 we hold all factors constant, except the design points (random design) and the noise process
 that were regenerated.\\
 E1 corresponds to the Bayesian estimator based on Gaussian
prior with large variance, E2 to the Bayesian estimator based on
Gaussian prior with small variance and E3 to the estimator built
following the Kerkyacharian and Picard procedure in \cite{Kerk2}.\\
 In order to implement E1, we made the following
choices of hyperparameters described in section 2.2.2 : in
(\ref{hyperparametrebigw}), $q_{1}=q_{2}=q=1$ proved to be a good
compromise whatever the function of interest to be estimated while
leading to good graphics reconstructions. We set $w(n)=20\times
n^{-q/2}$ and $\tau(n)=20\times\sigma^2/(n\log(n))$.
To implement E2, we set $c_{1}=1$, $c_{2}=2$, $\alpha=0.5$ and $\beta=1$, following the choices recommended in \cite{Abra2}.\\
The following plots compare the visual quality of the
reconstructions (see Fig. 3. to Fig. 8). The solid line is the
estimator and the dotted line is
the true function. \\

\begin{table}[htbp]
\begin{center}
\begin{tabular}[]{|c|c|c|c|c|c|c|c|}
\cline{3-8}
\multicolumn{2}{r|}{}&\multicolumn{3}{c|}{}&\multicolumn{3}{c|}{}\\
 \multicolumn{2}{r|}{}&\multicolumn{3}{c|}{RSNR=4}&\multicolumn{3}{c|}{RSNR=7}\\
 \multicolumn{2}{r|}{}&\multicolumn{3}{c|}{}&\multicolumn{3}{c|}{}\\
 \cline{2-8}
\multicolumn{1}{r|}{}&&&&&&&\\
\multicolumn{1}{r|}{}& design density& E1 & E2 & E3 & E1& E2 & E3\\
\multicolumn{1}{r|}{}&&&&&&&
\\
\hline \textit{Blocks}&
Sine & 0.0194 & 0.0219 &0.0227 &0.0113 & 0.0161&0.0129\\
&  Hole2 & 0.0196& 0.0220 &0.0226  &0.0114 &0.0163 &0.0130\\
\hline \textit{Bumps}&
Sine &0.0243 &0.240  &0.259  &0.0156 &0.0167 &0.0172\\
&Hole2 & 0.0241&0.0237  &0.0253 &0.0155 &0.0167 &0.0169\\
\hline \textit{HeaviSine}&
Sine & 0.0164& 0.0141& 0.0133& 0.0103&0.0092 &0.0093\\
&Hole2 &  0.0169& 0.0146 &0.0138  &0.0107 & 0.0097&0.0096 \\
\hline \textit{Doppler}&
Sine &  0.0236&0.0231 &0.0236  &0.0157 &0.0238 &0.0248\\
&Hole2 &0.0244  &0.0238  &0.0248  & 0.0166& 0.0172&0.0176  \\
 \hline
\end{tabular}
\end{center}
\caption{Values of RMSE over 100 runs} \label{Tableau}
\end{table}
\par
\centerline{\includegraphics[bb=60 190 545
590,scale=0.50]{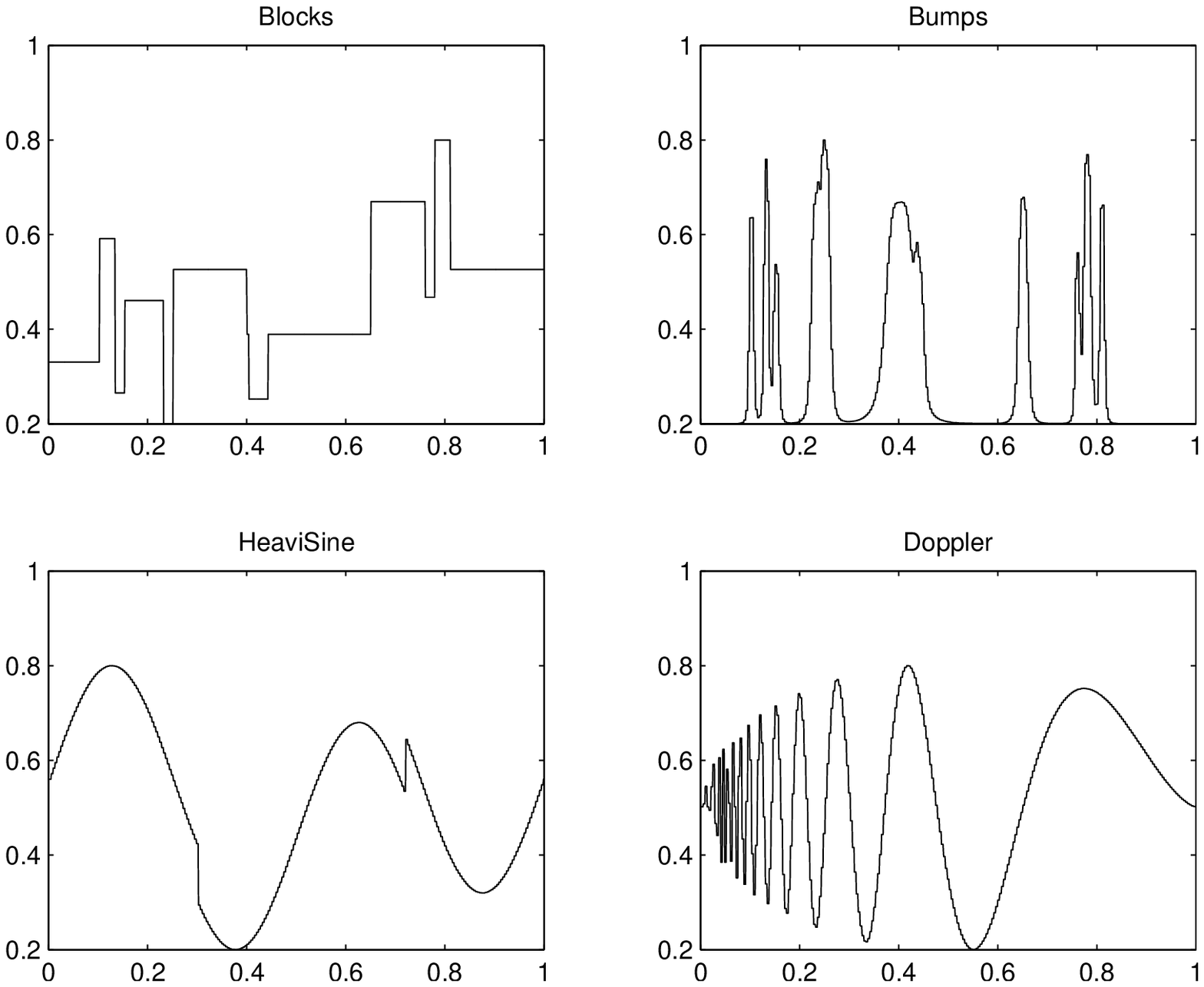}} \centerline{Fig. 1 Test functions}
\ \par \centerline{\includegraphics[bb=60 190 545
590,scale=0.50]{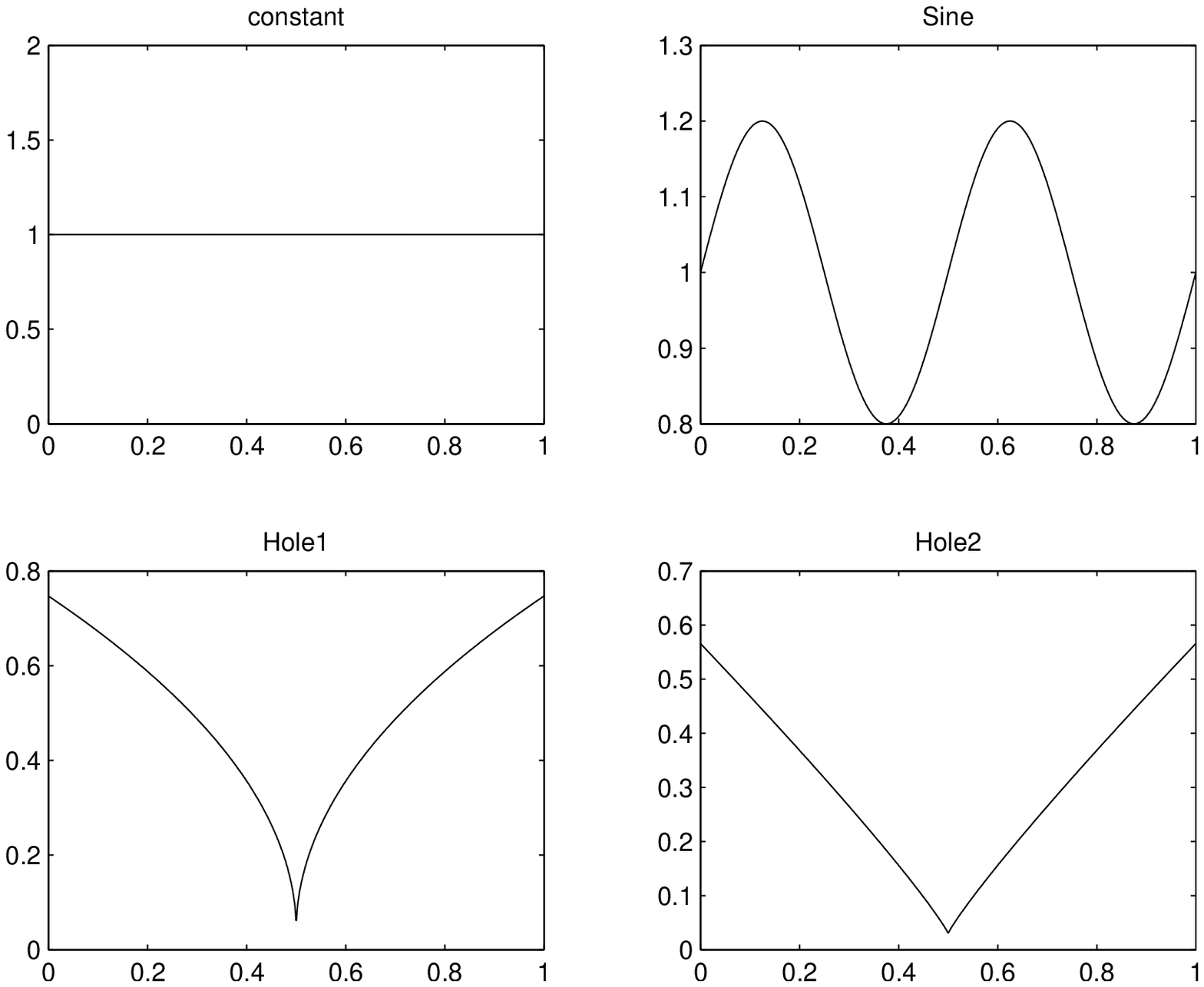}} \centerline{Fig. 2 Design
density}  \ \par \centerline{\includegraphics[bb=60 190 545
590,scale=0.50]{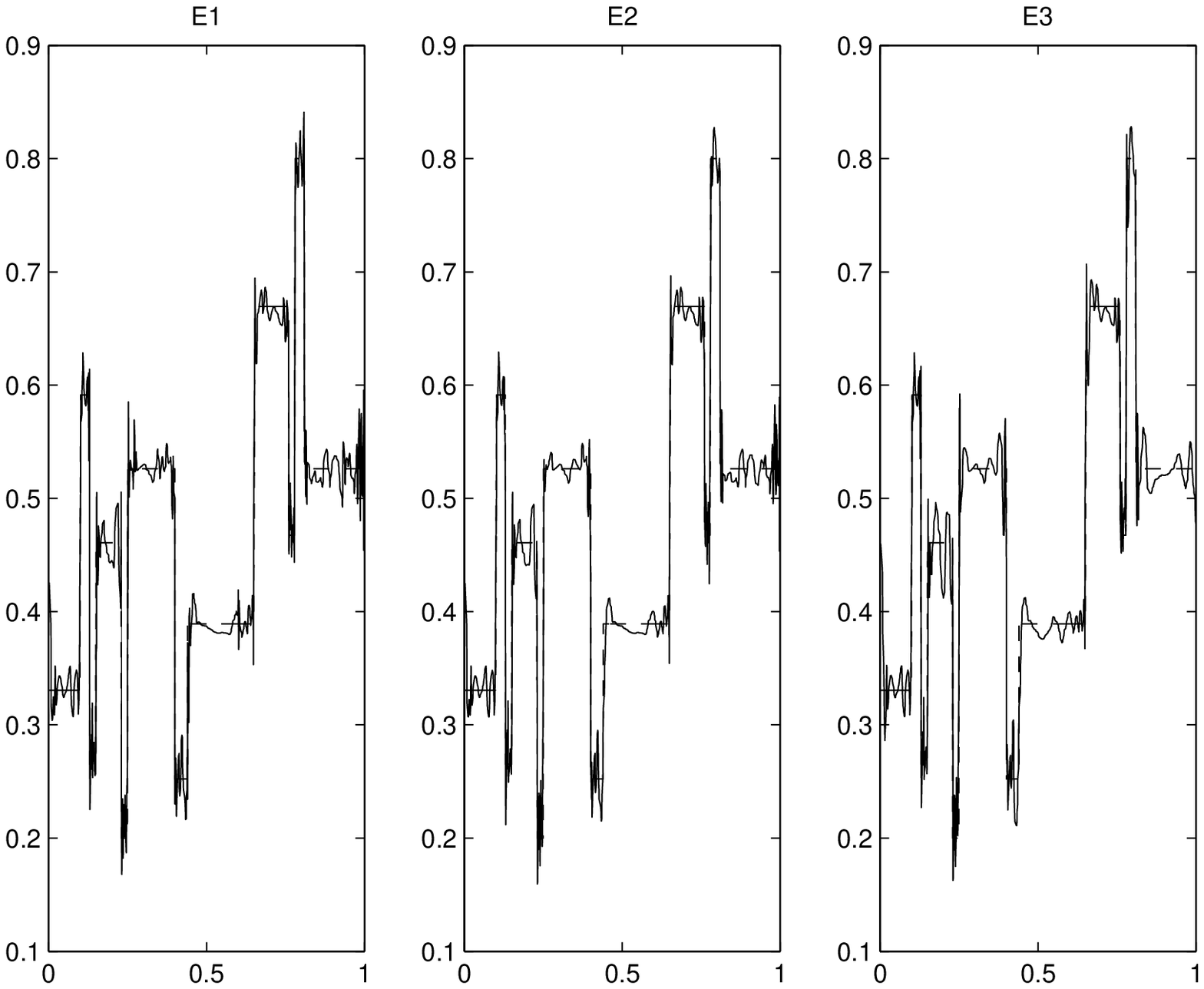}} \centerline{Fig. 3  Blocks target and
Sine density, RSNR=4} \ \par \centerline{\includegraphics[bb=60 190
545 590,scale=0.50]{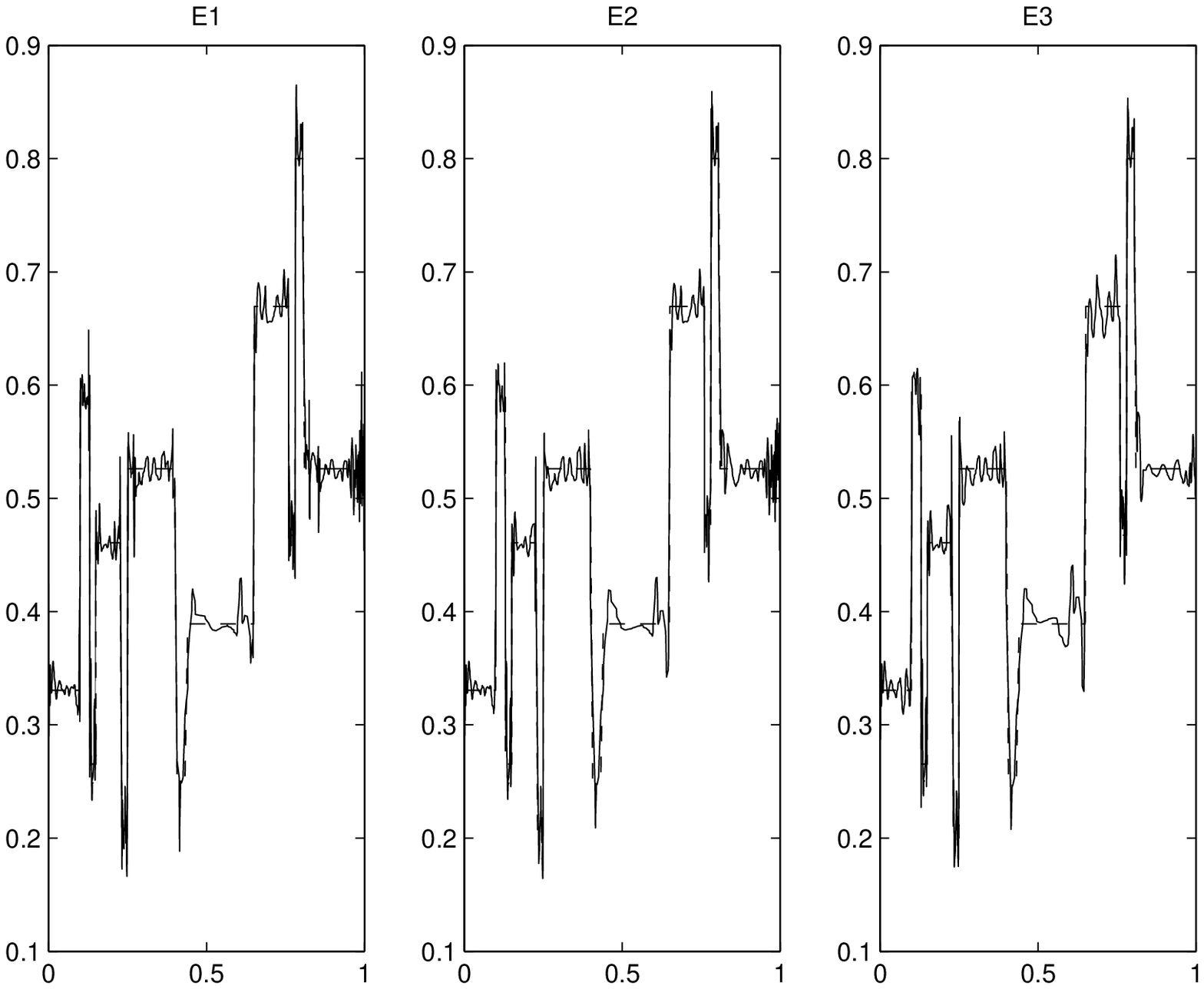}} \centerline{Fig. 4 Blocks target and
Hole2 design density, RSNR=4} \
\par
\centerline{\includegraphics[bb=60 190 545
590,scale=0.50]{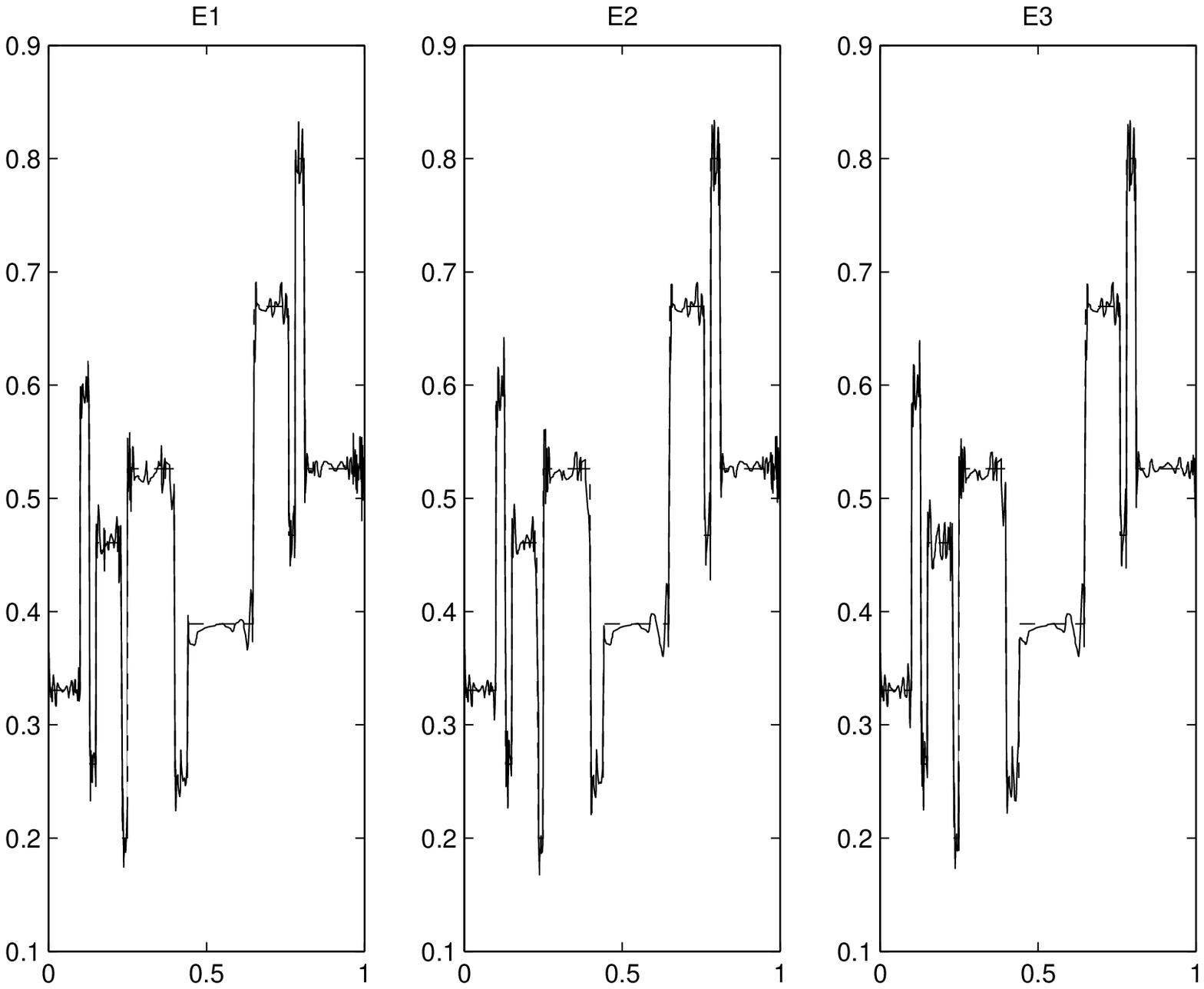}} \centerline{Fig. 5 Blocks target and
Hole2 design density,  RSNR=7} \
\par
\centerline{\includegraphics[bb=60 190 545
590,scale=0.50]{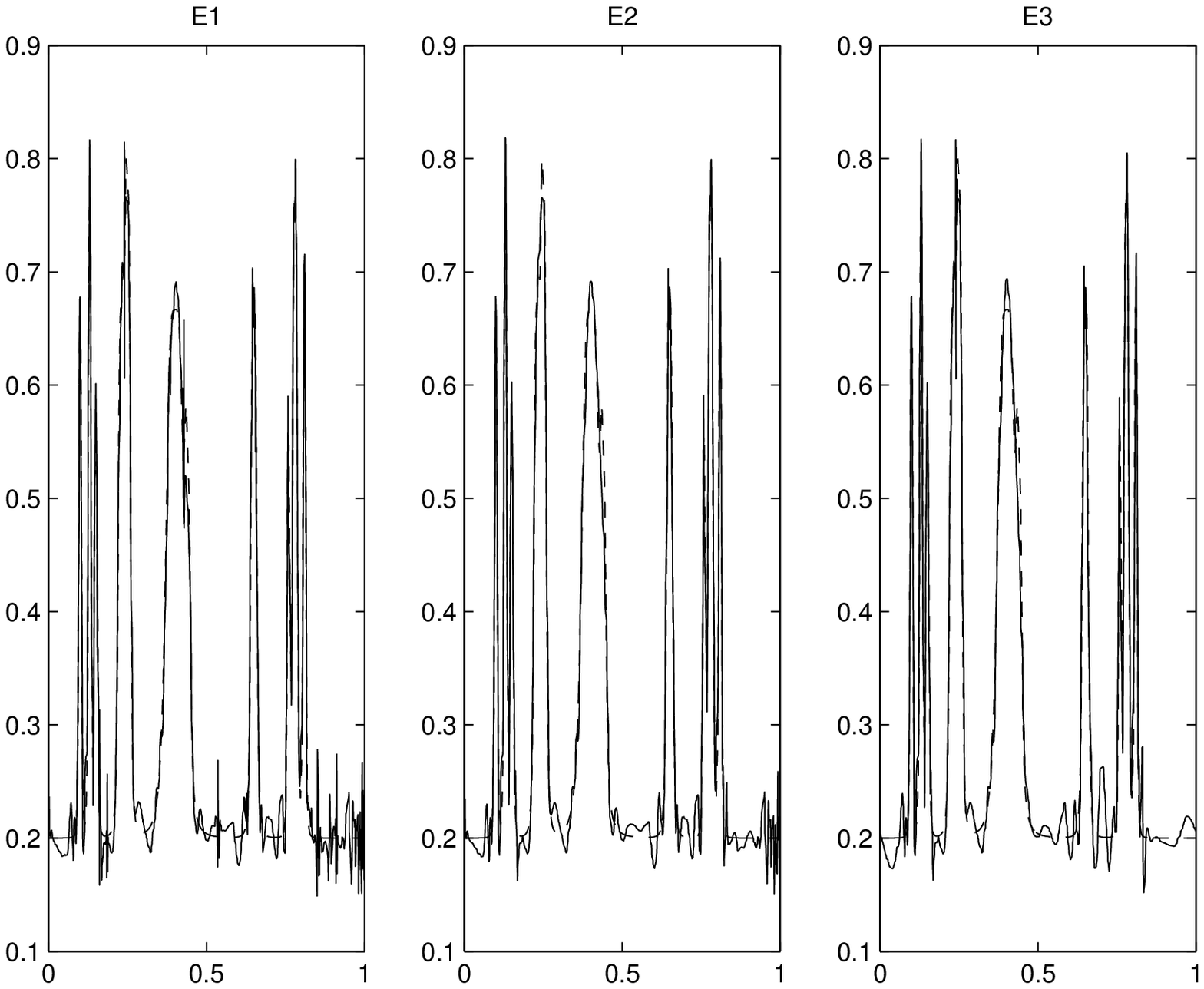}} \centerline{Fig. 6  Bumps target and
Sine design density,  RSNR=4} \
\par
\centerline{\includegraphics[bb=60 190 545
590,scale=0.50]{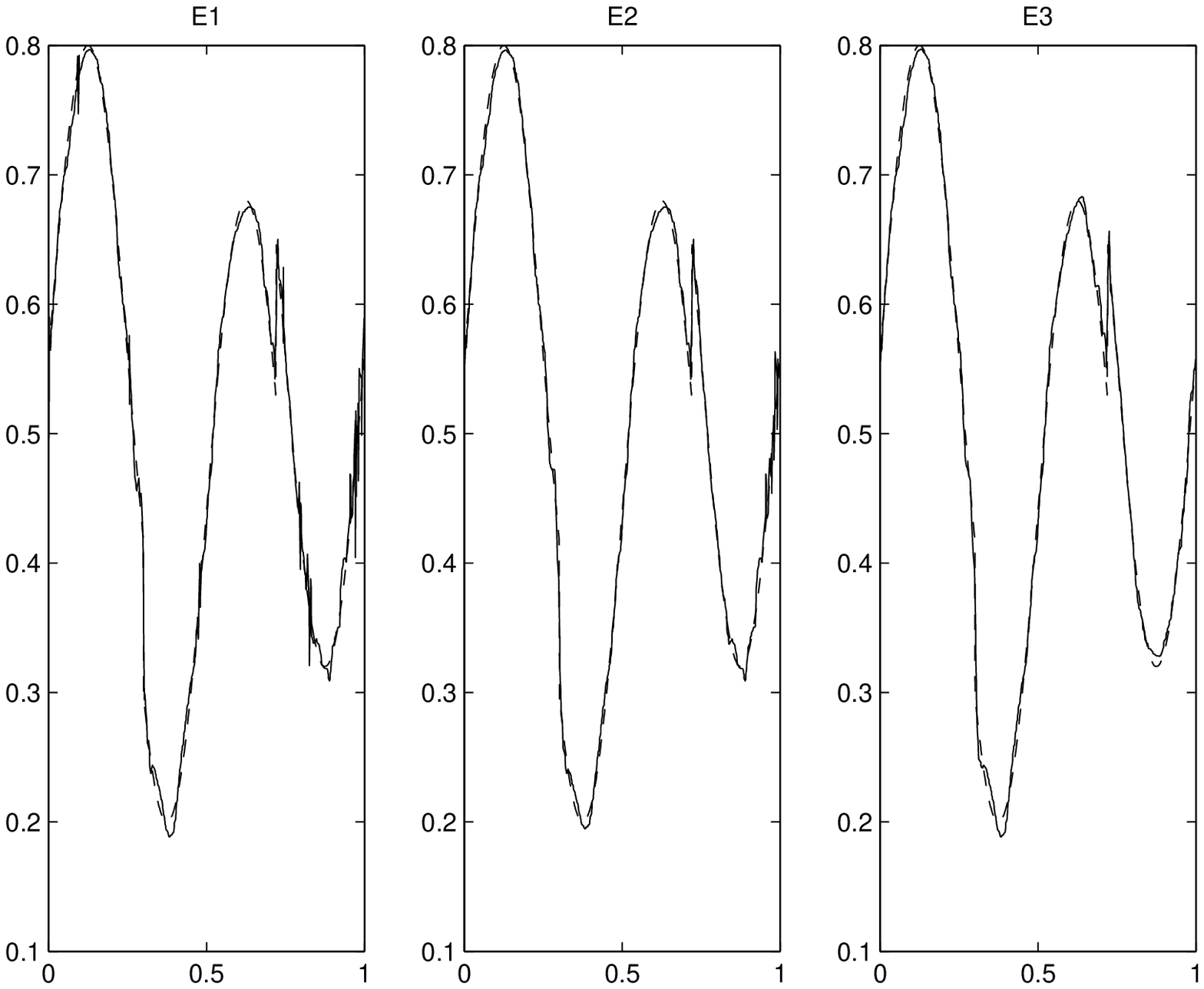}} \centerline{Fig. 7  HeaviSine target and
Sine design density, RSNR=7} \ \par
\centerline{\includegraphics[bb=60 190 545
590,scale=0.50]{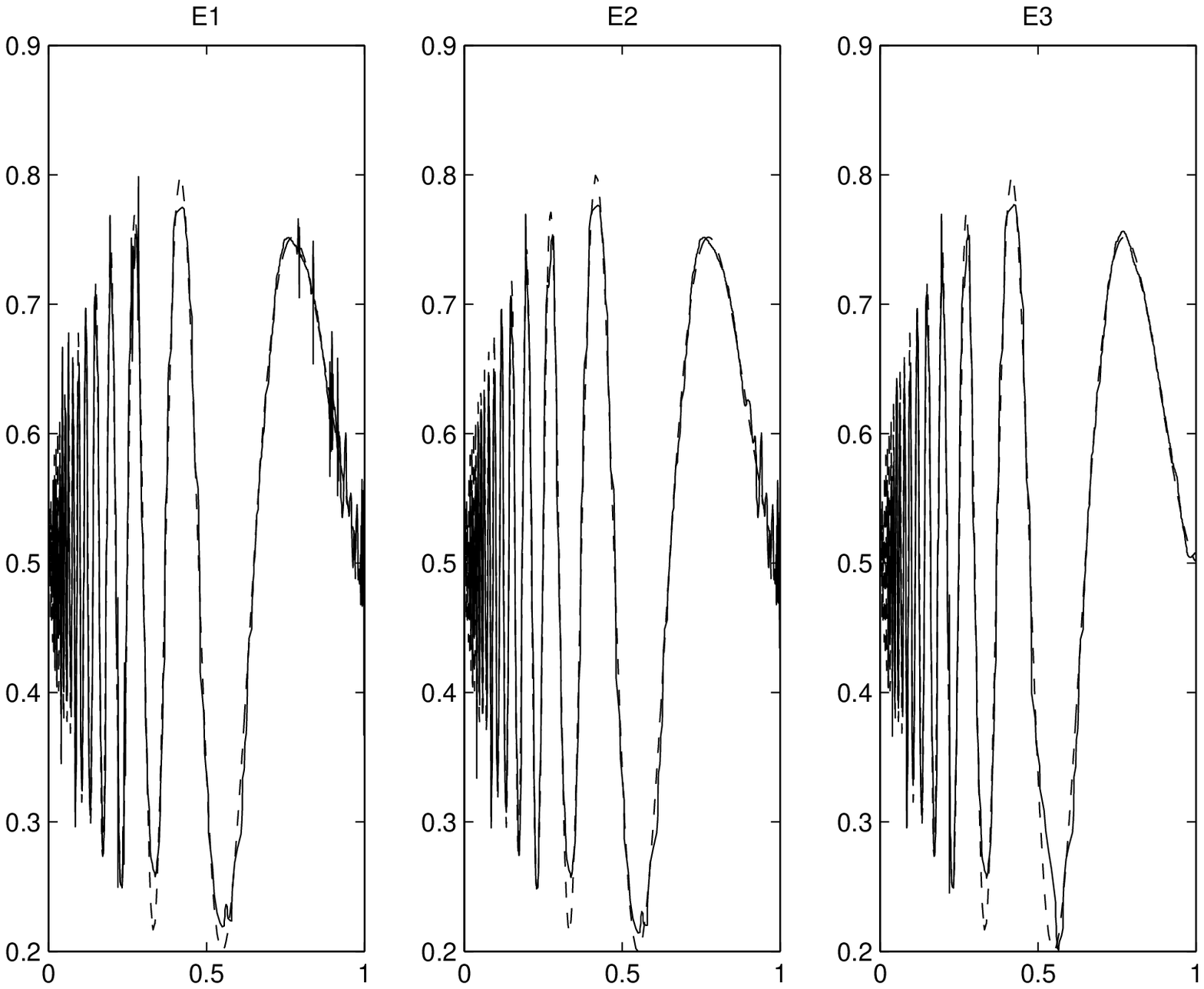}} \centerline{Fig. 8  Doppler target and
Hole2 design density, SNR=4} \
\par
\ \par We shall now comment and discuss the results displayed in
Tab.1 as well as the various visual reconstructions. \\ The
performances are always better for the Bayesian estimators except
for the case of the HeaviSine test function. More precisely, the
RMSE for Blocks whatever the noise level and design densities are
smaller for Estimator 1, moreover the RMSE are almost equal for
Estimator 1 and 2 in the case of Bumps test function, whatever the
design densities and for a noise level RSNR=4. This may be due to
the irregularity of the Bumps, Blocks and Doppler test functions
which are much rougher than the HeaviSine which is more regular.
Indeed, Estimator 1 and 2 tend to detect better the corner of
Blocks, the high peaks in Bumps, and the high frequency parts of
Doppler as the graphics show it. We may explain this by the fact
that Estimators 1 and 2 have level-dependent thresholds whereas
Estimator 3 has a hard universal threshold. \\ As for the
reconstructions, one can see that they are slighly better in the
case of Sine density and small noise, whereas there are small
deteriorations when a hole occurs in the design density but this
change does not affect the visual quality in too big proportions.
This fact highlights the interest of "warping" the wavelet basis.
Warping the basis allows the estimators to behave still correctly
when the design densities are far from the uniform density such as
in the case of Hole2.

\vspace{8mm} \noindent{\large\bf Acknowledgements}

The author wishes to thank her advisor Dominique Picard and Vincent
Rivoirard for interesting discussions and suggestions.
\\
\par

 \ \par \noindent
Laboratoire de Probabilit\'es et Mod\`eles al\'eatoires, UMR 7599,
Universit\'e Paris 6, case 188, 4, Pl. Jussieu, F-75252
Paris Cedex 5, France. \\
E-mail: thanh.pham\_ngoc@upmc.fr

\section{Appendix}
\vspace{6 mm} In the sequel $C$ denotes some positive constant which
may change from one line to another line. We also assume without
loss of generality that $\sigma=1$ in model (\ref{modele_de_base}).
\\We have that
$$\E(\psi_{jk}^{2}(G(X)))=\int\psi_{jk}^{2}(G(x))g(x)dx=\int \psi_{jk}^{2}(y)dy=1.$$
hence we get $\E(\gamma_{jk}^{2})=1/n$, the expression of
$\gamma_{jk}^{2}$ being given by (\ref{bruit_gamma}).\\
 Let us define the following event:
\begin{equation}\label{evenementOmega}
\Omega_{n}^{\delta}=\{|\gamma_{jk}^{2}-1/n| \leq \delta\}.
\end{equation}
To make proofs clearer we recall the Bernstein inequality that we
will use in the sequel. (see in \cite{Massart_concentration}
Proposition 2.8 and formula (2.16))
\begin{pro}\label{PropBernstein}
Let $Z_{1},\dots, Z_{n}$ be independant and square integrable random
variables such that for some nonnegative constant $b$, $Z_{i}\leq b$
almost surely for all $i \leq n$. Let
$$ S=\sum_{i=1}^{n} (Z_{i}-\E[Z_{i}])$$
and $v=\sum_{i=1}^{n} \E(Z_{i}^{2})$. Then for any positive $x$, we
have
$$ \p[S \geq x] \leq \exp\left( \frac{-v}{b^{2}}h(\frac{bx}{v}) \right)$$
where $h(u)=(1+u)\log(1+u)-u$. \\
It is easy to prove that
$$ h(u) \leq \frac{u^{2}}{2(1+u/3)}$$
which immediately yields
$$\p[S\geq x] \leq \exp\left( \frac{-x^{2}}{2(v+bx/3)} \right).$$
\end{pro}

\begin{lem}\label{lemma_0.}
Let $\varsigma$ be some positive constant. We have
\begin{equation}\label{Proba_gamma_J_alpha}
\mathbb{P}(|\gamma_{jk}^{2}-1/n|>\varsigma/n)\leq
2e^{-n^{1-1/\alpha}\frac{\varsigma^2}{2C(1+\varsigma/3)}}\quad
\forall\;j\leq J_{\alpha}
\end{equation}
\begin{equation}\label{Proba_gamma_J_n}
\mathbb{P}(|\gamma_{jk}^{2}-1/n|>\varsigma/n)\leq
2e^{-\varsigma^{2}\log(n)/(C\|\psi\|^{4}_{4}+\varsigma\|\psi\|^{2}_{\infty})}
\quad \forall\;j\leq J_{n}.
\end{equation}
\end{lem}

\vspace{6mm}\textbf{Proof of Lemma \ref{lemma_0.}} \\
Let us deal with the first case $j\leq J_{\alpha}$. To bound
$\p(|\gamma_{jk}^{2}-1/n|>\varsigma/n)$ we will use the Bernstein
inequality and apply Proposition \ref{PropBernstein}. In the present
situation
$Z_{i}=(1/n^2)\psi_{jk}^{2}(G(X_{i}))$. \\
First of all, in order to apply the Bernstein inequality, we need
the value of the sum
$$v=\sum_{i=1}^{n}\E[((1/n^{2})\psi_{j,k}^{2}(G(X_{i})))^{2}]$$
we have
\begin{eqnarray}\E\psi_{j,k}^{4}(G(X))&=& \int_{0}^{1}\psi_{j,k}^{4}(G(x))g(x)dx=\int_{0}^{1}\psi_{j,k}^{4}(y)dy \nonumber\\
&\leq& \int_{0}^{1}2^{2j}\psi^{4}(2^jy-k)dy\leq 2^j\int\psi^{4}(y)dy
\leq C\|\psi\|^{4}_{4}2^{j} \label{esperancePsi4}
\end{eqnarray}
hence
$$(1/n^{4})\sum_{i=1}^{n}\E\psi_{j,k}^{4}(G(X_{i}))\leq
(C/n^{3})2^{J_{\alpha}}=\frac{C}{n^{3-1/\alpha}}$$ \\
moreover
$$\psi_{jk}^{2}(G(X))\leq \|\psi\|^{2}_{\infty}2^{j}\leq Cn^{1/\alpha}\quad j\leq J_{\alpha}$$
so
\begin{equation*}
\p( |\gamma_{jk}^{2}-1/n)|> \varsigma/n) \leq
2\exp(-\frac{\varsigma^2}{2C(1+\varsigma/3)}\frac{n^{-2}}{n^{-3+1/\alpha}}).\end{equation*}
Let us now deal with the second case $j\leq J_{n}$. To bound
$\p(|\gamma_{jk}^{2}-1/n|>\varsigma/n)$ we will follow the lines of
the proof of the first case. Here again
$$ Z_{i}=1/n^2\psi_{jk}^{2}(G(X_{i})).$$
According to (\ref{esperancePsi4}), we have
$$\E(1/n^{4}\psi_{jk}^4(G(X)))\leq C2^j/n^4\leq C/(n^3\log(n)), \quad \; j\leq J_{n}$$
and $$ v=\sum_{i=1}^{n}\E(1/n^{4}\psi_{jk}^4(G(X)))\leq
C\|\psi\|^{4}_{4}/(n^2\log(n))$$ and
$$1/n^{2}\psi_{jk}^2(G(X))) \leq \|\psi\|^{2}_{\infty}2^{j}/(n^2) \leq \|\psi\|^{2}_{\infty}/(n\log(n)), \quad \; j\leq J_{n}$$
consequently
\begin{equation*} \p(|\gamma_{jk}^{2}-1/n|>\varsigma/n)\leq
2e^{-\varsigma^{2}\log(n)/(C\|\psi\|^{4}_{4}+\varsigma\|\psi\|^{2}_{\infty})}.\end{equation*}
\vspace{6mm} \\
The following lemma shows that the properties of the Bayesian
estimators $\check{f}$ and $\hat{f}$ can be controlled on the event
$\Omega_{n}^{\delta}$. To lighten the notations for the proof of
this lemma, we will denote $\Omega_{n}$ for $\Omega_{n}^{\delta}$
and $\Omega_{n}^{c}$ the complementary of $\Omega_{n}$.
\begin{lem}\label{lemma 3}
We have
$$\E[I(\Omega_{n}^{c})\|\check{f}-f\|^{2}_{2}]=o((\log(n)/n)^{2s/(2s+1)}$$
$$\E[I(\Omega_{n}^{c})\|\hat{f}-f\|^{2}_{2}]=o((1/n)^{1-1/\alpha}\log(n)).$$
\end{lem}

 \vspace{6mm}
\textbf{Proof of Lemma \ref{lemma 3}.} \\

We have
\begin{eqnarray*}
\E\left[I(\Omega_{n}^{c})\|\check{f}-f\|^{2}_{2}\right] &\leq&
CJ_{n}\E\left[\sum_{j\leq
J}\sum_{k}(\tilde{\beta}_{jk}-\beta_{jk})^{2}I(\Omega_{n}^{c})\right]
+\mathbb{P}(\Omega_{n}^{c})\sum_{j>J_{n}}\left(\sum_{k}\beta_{jk}^{2}\right)^{1/2}\\
&\leq& V+B.
\end{eqnarray*}

Let us first deal with the variance term V. The estimator
$\tilde{\beta}_{jk}$ can be written as
$\tilde{\beta}_{jk}=w_{jk}\hat{\beta}_{jk}$ with $0\leq w_{jk}\leq
1$. We have
\begin{eqnarray*}
V&\leq&C J_n\mathbb{E}\left[\sum_{j\leq J_n,k}\left(w_{jk}(\hat\beta_{jk}-\beta_{jk})-(1-w_{jk})\beta_{jk}\right)^2I({\Omega_{n}^c})\right]\\
&\leq&2C J_n\mathbb{E}\left[\sum_{j\leq J_n}\sum_kw_{jk}^2(\hat\beta_{jk}-\beta_{jk})^2I({\Omega_{n}^c})\right]+2C J_n\sum_{j\leq J_n}\sum_k\mathbb{E}\left[(1-w_{jk})^2\beta_{jk}^2I({\Omega_{n}^c})\right]\\
&\leq&2C J_n\mathbb{E}\left[\sum_{j\leq
J_n}\sum_k(\hat\beta_{jk}-\beta_{jk})^2I({\Omega_{n}^c})\right]+2C
J_n\sum_{j\leq
J_n}\sum_k\mathbb{E}\left[\beta_{jk}^2I({\Omega_{n}^c})\right]
\end{eqnarray*}
because $0\leq w_{jk}\leq 1$. Then, using Cauchy Scharwz inequality
we get
\begin{eqnarray*}
V&\leq&2C J_n\sum_{j\leq
J_n}\sum_k\left[\mathbb{E}(\hat\beta_{jk}-\beta_{jk})^4\right]^{\frac{1}{2}}\mathbb{P}(\Omega_{n}^c)^{\frac{1}{2}}+2C
J_n\sum_{j\leq J_n}\sum_k\beta_{jk}^2\mathbb{P}(\Omega_{n}^c).
\end{eqnarray*}
Using (\ref{Proba_gamma_J_n}) and (\ref{puissance_4_diff_beta_chap})
we have
\begin{eqnarray*}
V\leq
2CJ_{n}2^{J_{n}}e^{-\varsigma^{2}\log(n)/(2C\|\psi\|^{4}_{4}/n+\varsigma\|\psi\|^{2}_{\infty})}
+2CJ_{n}\|f(G^{-1})\|^{2}_{2}e^{-\varsigma^{2}\log(n)/(C\|\psi\|^{4}_{4}+\varsigma\|\psi\|^{2}_{\infty})}.
\end{eqnarray*}
We recall that $2^{J_n}=n/\log(n)$, accordingly by choosing
$\varsigma$ large enough we have
$$V=o((\log(n)/n)^{2s/(2s+1)}$$
As for the term $B$
$$B\leq
C2^{-2J_{n}s}e^{-\varsigma^{2}\log(n)/(C\|\psi\|^{4}_{4}+\varsigma\|\psi\|^{2}_{\infty})}$$
which completes the proof for $\check{f}$. \\
The proof for $\hat{f}$ is similar, all inequalities hold a fortiori
since, in the case of the estimator $\hat{f}$ we have
$\mathbb{P}(\Omega_{n}^{c})\leq e^{-Cn^{1-1/\alpha}}$ (see
(\ref{Proba_gamma_J_alpha})).

\vspace{6mm}

Let us place in the setting of Theorem 1. We recall that
$\tilde{\beta}_{jk}$ is zero whenever $|\hat{\beta}_{jk}|$ falls
below a threshold $\lambda_{B}$ and we have the following lemma
concerning the behavior of $\lambda_{B}$
\begin{lem}\label{lemma_1} On the event
$\Omega_{n}^{\delta}$ defined by (\ref{evenementOmega}) with
$\delta=1/(2n)$, for $\alpha>1$ we have
\begin{equation}\label{Lambda} \lambda_{B}\approx
\sqrt{\frac{\log(n)}{n}}, \quad j<J_{\alpha}\end{equation} and
$J_{\alpha}$ is taken such that
$2^{J_{\alpha}}=(\frac{3}{2n})^{-1/\alpha}$.
\end{lem}
\vspace{6mm}
\textbf{Proof of Lemma \ref{lemma_1}.} We follow the lines of the proof of lemma 1. in \cite{Abra1}.\\
On the one hand we have (see proof of lemma 1. in \cite{Abra1} page
228)
$${\lambda_{B}}^{2}\leq
\frac{2\gamma_{jk}^{2}(\gamma_{jk}^{2}+\tau_{j}^{2})}{\tau_{j}^{2}}\log\bigg(\frac{1-\pi_j}{\pi_j}\frac{\sqrt{\gamma_{jk}^{2}+\tau_{j}^2}}{\gamma_{jk}}+c\bigg)$$
where $c$ is some suitable positive constant. Besides, we have
$1/(2n) \leq \gamma_{jk}^{2} \leq 3/(2n)$, therefore
$${\lambda_{B}}^2 \leq \frac{2(3/(2n))((3/(2n))+c_{1}(3/(2n)))}{c_{1}(3/(2n))}\log\bigg(\frac{1-c_{2}(3/(2n))^{\beta/\alpha}}{c_{2}(3/(2n))^{\beta/\alpha}}
\frac{\sqrt{(1+c_{1})(3/(2n))}}{\sqrt{1/(2n)}}+c\bigg)$$ \\
hence we get
$$  {\lambda_{B}}^{2}\leq \tilde{c}(1/n)\log(\tilde{c}(1/n)^{(-\frac{\beta}{\alpha})}+c)$$
where $\tilde{c}$ denotes a positive constant depending on $c_{1}$
and $c_{2}$ and which may be different at different places. Since
$$\tilde{c}(1/n)\log(\tilde{c}(1/n)^{(-\frac{\beta}{\alpha})}+c)
\approx -\tilde{c}(\beta/\alpha)(1/n)\log(1/n)$$ we finally get
$${\lambda_{B}}^{2}\leq -\tilde{c}(\beta/\alpha)(1/n)\log(1/n).$$
On the other hand, for the reverse inequality, we have (see proof of
lemma 1. in \cite{Abra1} page 228 and formula (14) in \cite{Abra1}
page 221)
$${\lambda_{B}}^{2}\geq
\frac{2\gamma_{jk}^{2}(\gamma_{jk}^{2}+\tau_{j}^{2})}{\tau_{j}^{2}}\log\bigg(\frac{1-\pi_j}{\pi_j}\frac{\sqrt{\gamma_{jk}^{2}+\tau_{j}^2}}{\gamma_{jk}}\bigg)$$
but $|\gamma_{jk}^{2}-1/n|\leq\ 1/(2n)$ consequently one has
$${\lambda_{B}}^{2}\geq -\tilde{c}(\beta/\alpha)(1/n)(\log(1/n))$$
which completes the proof.

\vspace{6mm} \textbf{Proof of Theorem 1.} \\
Let us place on the event $\Omega_{n}^{\delta}$ defined by
(\ref{evenementOmega}) with $\delta=1/(2n)$. \\
 By the usual decomposition of the
MISE into a variance and a bias term we get
\begin{eqnarray} \E\|\hat{f}-f\|^{2}_{2}&\leq& 2\big[ \E\| \sum_{j \leq
J_{\alpha}}\sum_{k}(\tilde{\beta}_{jk}-\beta_{jk})\psi_{j,k}(G)\|^{2}_{2}
+
\|\sum_{j>J_{\alpha}}\sum_{k}\beta_{jk}\psi_{j,k}(G)\|^{2}_{2}\big] \nonumber\\
&\leq& 2(V+B) \nonumber
\end{eqnarray}
with
$$V=\E\| \sum_{j \leq
J_{\alpha}}\sum_{k}(\tilde{\beta}_{jk}-\beta_{jk})\psi_{j,k}(G)\|^{2}_{2}$$
$$B=\|\sum_{j>J_{\alpha}}\sum_{k}\beta_{jk}\psi_{j,k}(G)\|^{2}_{2}.$$

We first deal with the term $V$. We have
$$\| \sum_{j \leq
J_{\alpha}}\sum_{k}(\tilde{\beta}_{jk}-\beta_{jk})\psi_{j,k}(G)\|^{2}_{2}\leq
J_{\alpha}\sum_{j \leq J_{\alpha}}\|\sum_{k}
(\tilde{\beta}_{jk}-\beta_{jk})\psi_{j,k}(G)\|^{2}_{2}.$$ We want to
show that
$$\|\sum_{k}(\tilde{\beta}_{jk}-\beta_{jk})\psi_{j,k}(G)\|^{2}_{2}\leq C\sum_{k}(\tilde{\beta}_{jk}-\beta_{jk})^{2}.$$
For this purpose we have
\begin{eqnarray*}
\| \sum_{k}(\tilde{\beta}_{jk}-\beta_{jk})\psi_{j,k}(G)\|^{2}_{2}
&=& \int|\sum_{k}
(\tilde{\beta}_{jk}-\beta_{jk})\psi_{jk}(G(x))|^2 dx\\
&=&\int|\ \sum_{k} (\tilde{\beta}_{jk}-\beta_{jk})\psi_{jk}(x)|^2
\frac{1}{g(G^{-{1}}(x))}dx\\
&=& \|
\sum_{k}(\tilde{\beta}_{jk}-\beta_{jk})\psi_{j,k}\|^{2}_{\mathbb{L}_2(\omega
)}
\end{eqnarray*}
where $\omega(x)=1/(g(G^{-1}))(x)$.\\ Now using inequality $(44)$ p.
$1075$ in \cite{Kerk2} we have
$$ \|\sum_{k}(\tilde{\beta}_{jk}-\beta_{jk})\psi_{j,k}\|^{2}_{\mathbb{L}_2(\omega
)}\leq
C2^{j}\sum_{k}|\tilde{\beta}_{jk}-\beta_{jk}|^{2}\omega(I_{j,k})$$
where $I_{j,k}$ denotes the interval
$[\frac{k}{2^j},\frac{k+1}{2^j}]$ and
$\omega(I_{jk})=\int_{I_{jk}}\omega(x)dx$. But the design density
$g$ is bounded below by $m$. Hence we get
$$ \omega(I_{j,k})\leq 2^{-j}/m$$ and consequently
$$ \|
\sum_{k}(\tilde{\beta}_{jk}-\beta_{jk})\psi_{j,k}\|^{2}_{\mathbb{L}_2(\omega
)} \leq C\sum_{k}(\tilde{\beta_{jk}}-\beta_{jk})^{2}.$$

We decompose now $V$ into three terms
$$ V\leq C J_{\alpha}\E\sum_{j\leq J_{\alpha}}\sum_{k}
[(\tilde{\beta}_{jk}-\beta_{jk}^{'})^{2}+(\beta_{jk}^{'}-\beta_{jk}^{''})^{2}+(\beta_{jk}^{''}-\beta_{jk})^{2}]$$
where
$$\beta_{jk}^{'}=b_{j}\hat{\beta}_{jk}I\{|\hat{\beta}_{jk}|\geq
\kappa\lambda_{B}\}$$ with $\kappa$ a positive constant and
$$ b_{j}=\frac{\tau_{j}^{2}}{\tau_{j}^{2}+\gamma_{jk}^{2}}$$
$$\beta_{jk}^{''}=b_{j}\beta_{jk}.$$
As a consequence we have
$$V\leq CJ_{\alpha}(A_{1}+A_{2}+A_{3}).$$
We are now going to upperbound each term $A_{1}$, $A_{2}$ and
$A_{3}$. We start with $A_{1}$
$$A_{1}=\sum_{j\leq
J_{\alpha}}\sum_{k}\E[(\tilde{\beta}_{jk}-\beta^{'}_{jk})^{2}].$$ \\
As precised in the beginning of section 2.2 p 6,
$\tilde{\beta}_{jk}=0$ for $|\hat{\beta}_{jk}|<\lambda_{B}$. As
well, $\beta_{jk}^{'}=0$ for $|\hat{\beta}_{jk}|<\kappa \lambda_{B}$
and $\tilde{\beta}_{jk}-\beta_{jk}^{'}\rightarrow 0$ monotonically
as $\hat{\beta}_{jk} \rightarrow \infty$. Hence
$$\max_{\hat{\beta}_{jk}}|\tilde{\beta}_{jk}-\beta_{jk}^{'}|=b_{j}\lambda_{B}$$
which implies
$$A_{1}\leq C\sum_{j<J_{\alpha}}\sum_{k}\E(b_{j}^{2}\lambda_{B}^{2}). $$\\
We have $\lambda_{B}\approx \sqrt{\frac{\log n }{n}}$ and $b_{j}\leq
1$ for $j \leq J_{\alpha}$ hence we get
$$A_{1} \leq C\sum_{j\leq J_{\alpha}} \sum_{k=0}^{2^j-1}
\frac{\log(n)}{n}$$ \\
so
\begin{eqnarray} A_{1} &\leq& C\frac{\log(n)}{n} \sum_{j\leq
J_{\alpha}}2^j
\\ &\leq&
C\frac{\log(n)}{n}\bigg({\frac{1}{n}}\bigg)^{-1/\alpha}
\end{eqnarray}
finally
$$A_{1}=\mathcal{O}\big(\log(n)(\frac{1}{n})^{1-1/\alpha}\big)$$
Let us now consider the second term $A_{2}$
\begin{eqnarray} A_{2}&=& \sum_{j\leq J_{\alpha}} \sum_{k=0}^{2^j-1}
\E(\beta_{jk}^{'}-\beta_{jk}^{''})^{2} \nonumber \\
&=& \sum_{j\leq J_{\alpha}} \sum_{k=0}^{2^j-1}
\E(b_{j}\hat{\beta}_{jk}I\{|\hat{\beta}_{jk}|\geq \kappa\lambda_{B}
\}-b_{j}\beta_{jk})^{2} \nonumber \end{eqnarray} \\
We have that $b_{j} \leq 1$, consequently it follows
 \begin{eqnarray}
A_{2} &=& \sum_{j\leq J_{\alpha}} \sum_{k=0}^{2^j-1}
\E((\hat{\beta}_{jk}-\beta_{jk})^{2}I\{|\hat{\beta}_{jk}|\geq
\kappa\lambda_{B}\})+\E\sum_{j\leq J_{\alpha}} \sum_{k=0}^{2^j-1}
\beta_{jk}^{2}I\{|\hat{\beta}_{jk}|<\kappa\lambda_{B} \} \nonumber
\\&=& A_{2}^{'}+A_{2}^{''} \nonumber
\end{eqnarray}
We have
$$A_{2}^{'} \leq \sum_{j\leq J_{\alpha}} \sum_{k=0}^{2^j-1}
\E(\hat{\beta}_{jk}-\beta_{jk})^{2}$$ Using inequality $(64)$ in
\cite{Kerk2} p. 1086 we have
\begin{equation}\label{Inegalite_diff_carre_bernoulli}
\E(\hat{\beta}_{jk}-\beta_{jk})^{2} \leq
C\frac{1+\|f\|^{2}_{\infty}}{n}\end{equation} hence
$$A_{2}^{'}=\mathcal{O}((1/n)^{1-1/\alpha}).$$
We now bound the term $A_{2}^{''}$.
\begin{eqnarray}
A_{2}^{''}&=& \E\sum_{j\leq J_{\alpha}} \sum_{k=0}^{2^j-1}
\beta_{jk}^{2}I\{|\hat{\beta}_{jk}|<\kappa\lambda_{B}\}(I\{|\beta_{jk}|<2\kappa\lambda_{B}\}+I\{|\beta_{jk}|>2\kappa\lambda_{B}\}) \nonumber \\
&\leq& \E\sum_{j\leq J_{\alpha}}
\sum_{k=0}^{2^j-1}\beta_{jk}^{2}I\{|\beta_{jk}|<2\kappa\lambda_{B}\}+
\sum_{j\leq J_{\alpha}}
\sum_{k=0}^{2^j-1}\beta_{jk}^{2}\p(|\hat{\beta}_{jk}-\beta_{jk}|>\kappa\lambda_{B}) \nonumber\\
&=& T_{3}+T_{4}
\end{eqnarray}
We have
\begin{eqnarray}
T_{3}\leq C\sum_{j\leq J_{\alpha}}\lambda_{B}^{2}2^{j}\leq
C\frac{\log(n)}{n}{n}^{1/\alpha}=C\log(n)n^{-1+1/\alpha}. \nonumber
\end{eqnarray}
Let us focus on $T_4$, we have
$$\hat{\beta}_{jk}-\beta_{jk}=1/n\sum_{i=1}^{n}\psi_{j,k}(G(X_{i}))(f(X_{i})+\varepsilon_{i})-\mathbb{E}\psi_{j,k}(G(X))f(X)$$
Hence
\begin{equation*} \mathbb{P}(|\hat{\beta}_{jk}-\beta_{jk}|>
\kappa\sqrt{\log(n)/n})\leq
\mathbb{P}_{1}+\mathbb{P}_{2}\end{equation*} where
\begin{equation}\label{P_1}\mathbb{P}_{1}=\mathbb{P}(|1/n\sum_{i=1}^{n}\psi_{j,k}(G(X_{i}))(f(X_{i}))-\mathbb{E}\psi_{j,k}(G(X))f(X)|>
\kappa/2\sqrt{(\log(n)/n)})\end{equation} and
\begin{equation}\label{P_2}\mathbb{P}_{2}=\mathbb{P}(|1/n\sum_{i=1}^{n}\psi_{j,k}(G(X_i))\varepsilon_i|>\kappa/2\sqrt{(\log(n)/n)})\end{equation}
Kerkyacharian and Picard in \cite{Kerk2} in order to prove
inequality $(65)$ in \cite{Kerk2} showed p. 1088 that
\begin{equation}\label{P_1_valeur}\mathbb{P}_{1} \leq
2\exp(-\frac{3\kappa^{2}\log(n)}{4\|f\|_{\infty}(3+\kappa)})\end{equation}
if $2^j\leq n/\log(n)$. As for $\mathbb{P}_{2}$, conditionally on
$(X_{1},\dots,X_{n})$ we have
$$1/n\sum_{i=1}^{n}\psi_{j,k}(G(X_i))\varepsilon_i \sim N(0,\gamma_{jk}^{2})$$
where $\gamma_{jk}^{2}$ has been defined in (\ref{bruit_gamma}).
Using exponential inequality for Gaussian random variable we have
\begin{eqnarray} \mathbb{P}_{2}&\leq& \E(\exp(-\frac{\kappa^{2}\log(n)}{8n\gamma_{jk}^2})) \nonumber
\\&=&
\E
e^{-\frac{\kappa^{2}\log(n)}{8n\gamma_{jk}^2}}(I(|\gamma_{jk}^{2}-1/n|\leq
1/2n)+I(|\gamma_{jk}^{2}-1/n|>1/(2n))) \nonumber
\\ &\leq&
e^{-\frac{\kappa^{2}\log(n)}{12}}+\p(|\gamma_{jk}^{2}-1/n|>1/(2n)).
\label{P_2_valeur}
\end{eqnarray}
Using (\ref{Proba_gamma_J_alpha}) with $\varsigma=1/2$, we have for
$\alpha>1$
\begin{eqnarray}
T_{4}&\leq&
(2e^{(-Cn^{1-1/\alpha})}+e^{-\frac{\kappa^{2}\log(n)}{12}}+2\exp(\frac{-3\kappa^2\log(n)}{4\|f\|_{\infty}(3+\kappa)})
)\sum_{j\leq J_{\alpha}}
\sum_{k=0}^{2^j-1} \beta_{jk}^{2} \nonumber\\
&\leq& (2e^{(-Cn^{1-1/\alpha})}+ e^{-\frac{\kappa^{2}\log(n)}{12}}+
2\exp(\frac{-3\kappa^2\log(n)}{4\|f\|_{\infty}(3+\kappa)}))\|f(G^{-1})\|^{2}_{2}
\nonumber
\end{eqnarray}
It remains to fix $\kappa$ large enough so that we get
$$T_{4}= \mathcal{O}(\log(n)n^{-1+1/\alpha}).$$
So we have for $A_{2}^{''}$, with $\alpha>1$,
\begin{eqnarray*}
A_{2}^{''} = \mathcal{O} (\frac{\log(n)}{n^{1-1/\alpha}})
\end{eqnarray*}
Finally we get for $A_{2}$
$$A_{2}=\mathcal{O}(\log(n)(\frac{1}{n})^{1-1/\alpha}).$$
Let us consider now the term $A_{3}$
\begin{eqnarray}
A_{3}&\leq& C\sum_{j\leq J_{\alpha}} \sum_{k=0}^{2^j-1}
\E(\beta_{jk}^{"}-\beta_{jk})^{2} \nonumber \\
&=& C\sum_{j\leq J_{\alpha}}
\sum_{k=0}^{2^j-1}\beta_{jk}^{2}(1-b_{j})^{2}=\sum_{j\leq
J_{\alpha}}\sum_{k=0}^{2^j-1}(\frac{\gamma_{jk}^{2}}{\tau_{j}^{2}+\gamma_{jk}^{2}})^{2}\beta_{jk}^{2}.
\nonumber
\end{eqnarray}
Since $|\gamma_{jk}^{2}-1/n|\leq 1/(2n)$, we get
\begin{eqnarray}
A_{3}\leq \sum_{j\leq J_{\alpha}}
\left(\frac{3/(2n)}{c_{1}2^{-j\alpha}+1/(2n)}\right)^{2}\sum_{k=0}^{2^j-1}\beta_{jk}^{2}
\nonumber
\end{eqnarray}
but $f(G^{-1})$ belongs to the Besov ball $B^{s}_{2,\infty}(R)$
which entails
$$\sum_{k=0}^{2^j-1}\beta_{jk}^{2}\leq M2^{-2js}$$
hence
\begin{eqnarray}
A_{3}\leq C/n^{2}\sum_{j\leq J_{\alpha}}2^{2j(-s+\alpha)} \nonumber
\end{eqnarray}
We have
\begin{eqnarray}
A_{3}\leq
C/n^{2}(1/n)^{\frac{-2(-s+\alpha)}{\alpha}}=\mathcal{O}{(1/n)^{2s/\alpha}}.\nonumber
\end{eqnarray}
We are now in position to give an upper bound for the variance term
$V$ namely
\begin{equation}\label{Variance}
V \leq CJ_{\alpha}(\log(n)(1/n)^{1-1/\alpha}+(1/n)^{2s/\alpha}).
 \nonumber
\end{equation}
It remains to bound the bias term $B$. In \cite{Kerk2} p.1083 using
inequality (44) the authors have proved that for any $l$ we get
\begin{equation}\label{kerk_bernoulli_2}\|\sum_{j\geq
l}\sum_{k}\beta_{jk}\psi_{j,k}(G)\|_{2}\leq \sum_{j\geq
l}\|\sum_{k}\beta_{jk}\psi_{j,k}(G)\|_{2}\ \leq C\sum_{j\geq
l}2^{j/2}\bigg(\sum_{k}|\beta_{jk}|^{2}\omega(I_{j,k})\bigg)^{1/2}.\end{equation}
Applying (\ref{kerk_bernoulli_2}) with in our case of lower bounded
design density, $\omega(I_{j,k})\leq 2^{-j}/m$ and $l=J_{\alpha}$,
it follows
\begin{eqnarray*}
\|\sum_{j\geq J_{\alpha}}\sum_{k}\beta_{jk}\psi_{j,k}(G)\|_{2}
&\leq& C \sum_{j\geq J_{\alpha}}\bigg(\sum_{k}|\beta_{jk}|^{2}
\bigg)^{1/2}
\\
&\leq& C\sum_{j\geq J_{\alpha}}2^{-js} \leq C2^{-J_{\alpha}s}
\end{eqnarray*}
hence
$$B=\|\sum_{j>J_{\alpha}}\sum_{k=0}^{2^j-1}\beta_{jk}\psi_{j,k}(G(x))\|^{2}_{2} \leq C2^{-2J_{\alpha}s}
=C(1/n)^{2s/\alpha} $$ which completes the proof of Theorem 1.

\begin{lem}\label{lemma 2.}
Let $w_{jk}$ a sequence of random weights lying in $[0,1]$. We
assume that there exist positive constants $c$, $m$ and $K$ such
that for any $\eps>0$,
$$\check{\beta}_{n}=(w_{jk}\hat{\beta}_{jk})_{jk}$$
is a shrinkage rule verifying for any $n$,
\begin{equation}
w_{jk}(n)=0, \quad a.e. \quad \forall\, j\geq J_{n} \;
\textrm{with}\; 2^{J_{n}}\sim n/\log(n):={t}_n^{2}, \, \forall \,k
\end{equation}
\begin{equation} \label{lemme2condition2}
|\hat{\beta}_{jk}|\leq mt_{n} \Rightarrow w_{jk}\leq ct_{n},\quad
\forall \,j\leq J_{n},\, \forall \,k,
\end{equation}
\begin{equation} \label{lemme2condition3}
(1-w_{jk}(n))\leq K(\frac{t_{n}}{|\hat{\beta}_{jk}|}+t_{n}) \quad
a.e. \quad \forall \,j\leq J_{n},\, \forall \,k.
\end{equation}
and let
$$\check{f}=\sum_{j<J_{n}}\sum_{k}
w_{jk}\hat{\beta}_{jk}\psi_{jk}(G(x))$$ Then
$$\sup_{f(G^{-1})\in B^{s}_{2,\infty}(R)}\E\|\check{f}-f\|^{2}_{2}\leq
(\log(n)/n)^{2s/(2s+1)}.$$
\end{lem}

\vspace{6mm} \textbf{Proof of Lemma \ref{lemma 2.}.}
\begin{eqnarray}
\E\|\check{f}-f\|^{2}_{2}&\leq& 2C (J_{n}\sum_{j\leq
J_{n}}\sum_{k}\E(\check{\beta}_{jk}-\beta_{jk})^2+
\|\sum_{j>J_{n}}\sum_{k}\beta_{jk}^{2}\psi_{jk}(G(x))\|^{2}_{2})
\nonumber
\\ &\leq&
V_{1}+B_{1}.  \nonumber
\end{eqnarray}
We first consider the term $V_{1}$
\begin{eqnarray*}
V_{1}&\leq & 2J_n\E\sum_{j\leq
J_{n}}\sum_{k}(w_{jk}^{2}(\hat{\beta}_{jk}-\beta_{jk})^{2}+(1-w_{jk})^{2}\beta_{jk}^{2})I\{|\hat{\beta}_{jk}|\leq
mt_{n}\} \\&+& J_n\E\sum_{j\leq
J_{n}}\sum_{k}(w_{jk}^{2}(\hat{\beta}_{jk}-\beta_{jk})^{2}+(1-w_{jk})^{2}\beta_{jk}^{2})I\{|\hat{\beta}_{jk}|>
mt_{n}\}\\ &=& V_{1}^{'}+V_{1}^{"}
\end{eqnarray*}
$$V_{1}^{'}=J_{n}(T_{5}+T_{6})$$
$$T_{5}=\E\sum_{j\leq
J_{n}}\sum_{k}w_{jk}^{2}(\hat{\beta}_{jk}-\beta_{jk})^{2}I\{|\hat{\beta}_{jk}|\leq
mt_{n}\}$$ but according to (\ref{Inegalite_diff_carre_bernoulli})
we have for $2^j\leq \log(n)/n$
$$\E(\hat{\beta}_{jk}-\beta_{jk})^{2}\leq C\frac{1+\|f\|^{2}_{\infty}}{n} $$ hence
using (\ref{lemme2condition2}) it follows
$$T_{5} \leq Ct_{n}^{2}2^{J_{n}}1/n .$$
As for $T_{6}$
\begin{eqnarray}
T_{6}&=&\E\sum_{j\leq J_{n}}\sum_{k}
(1-w_{jk})^{2}\beta_{jk}^{2}I\{|\hat{\beta}_{jk}|\leq mt_{n}\} \nonumber \\
&\leq& \E\sum_{j\leq J_{n}}\sum_{k}
(1-w_{jk})^{2}\beta_{jk}^{2}I\{|\hat{\beta}_{jk}|\leq mt_{n}\}[I\{
 |\beta_{jk}|\leq 2mt_{n}\}+I\{|\beta_{jk}>2mt_{n}|\}]. \nonumber
 \end{eqnarray}
 By (\ref{weakbesov}) we get
\begin{eqnarray*}
 T_{6}\leq& 2(mt_{n})^{2s/(2s+1)}\|f\|^{2}_{W_{2/(1+2s)}} + \sum_{j\leq J_{n}}\sum_{k}\beta_{jk}^{2}\p(|\hat{\beta}_{jk}-\beta_{jk}|>
 mt_{n}).
\end{eqnarray*}
We are going to bound $\p(|\hat{\beta}_{jk}-\beta_{jk}|> mt_{n})$.
We have
\begin{equation*} \mathbb{P}(|\hat{\beta}_{jk}-\beta_{jk}|\geq
m\sqrt{\log(n)/n})\leq \mathbb{P}_{3}+\mathbb{P}_{4}\end{equation*}
where
\begin{equation}\label{P_3}\mathbb{P}_{3}=\mathbb{P}(|1/n\sum_{i=1}^{n}\psi_{j,k}(G(X_{i}))(f(X_{i})-\mathbb{E}\psi_{j,k}(G(X))f(X)|\geq
m/2\sqrt{\log(n)/n})\end{equation} and
\begin{equation}\label{P_4}\mathbb{P}_{4}=\mathbb{P}(|1/n\sum_{i=1}^{n}\psi_{j,k}(G(X_i))\varepsilon_i|>m/2\sqrt{(\log(n)/n)}).\end{equation}
Kerkyacharian and Picard in \cite{Kerk2} in order to prove
inequality $(65)$ in \cite{Kerk2} showed p. 1088 that
\begin{equation}\label{P_3_valeur}\mathbb{P}_{3} \leq
2\exp(-\frac{3m^{2}\log(n)}{4\|f\|_{\infty}(3+m)})\end{equation} if
$2^j\leq n/\log(n)$. As for $\mathbb{P}_{4}$, conditionally on
$(X_{1},\dots,X_{n})$ we have
$$1/n\sum_{i=1}^{n}\psi_{j,k}(G(X_i))\varepsilon_i \sim N(0,\gamma_{jk}^{2})$$
where $\gamma_{jk}^{2}$ has been defined in (\ref{bruit_gamma}).
\begin{eqnarray} \mathbb{P}_{4}&\leq& \E(\exp(-\frac{m^{2}\log(n)}{8n\gamma_{jk}^2})) \nonumber
\\&=&
\E
e^{-\frac{m^{2}\log(n)}{8n\gamma_{jk}^2}))}(I(|\gamma_{jk}^{2}-1/n|\leq
\varsigma/n)+I(|\gamma_{jk}^{2}-1/n|>\varsigma/n)) \nonumber
\\ &\leq&
e^{-\frac{m^{2}\log(n)}{8(\varsigma+1)}}+\p(|\gamma_{jk}^{2}-1/n|>\varsigma/n).
\label{P_4_valeur}
\end{eqnarray}
Using (\ref{Proba_gamma_J_n}) to bound
$\p(|\gamma_{jk}^{2}-1/n|>\varsigma/n)$ we get
\begin{eqnarray*}
\p(|\hat{\beta}_{jk}-{\beta}_{jk}|> mt_{n})&\leq&
2e^{-\varsigma^{2}\log(n)/(C\|\psi\|^{4}_{4}+\varsigma\|\psi\|^{2}_{\infty})}+e^{-\frac{m^{2}\log(n)}{8(\varsigma+1)}}+2e^{(-\frac{3m^{2}\log(n)}{4\|f\|_{\infty}(3+m)}}
\label{bernstein2}
\end{eqnarray*}
thus
\begin{equation}\label{Proba_Beta_jk_chapeau}
\p(|\hat{\beta}_{jk}-{\beta}_{jk}|> mt_{n})\leq
2n^{\frac{-\varsigma^2}{C\|\psi\|^{4}_{4}+\varsigma\|\psi\|^{2}_{\infty}}}+n^{\frac{-m^2}{8(\varsigma+1)}}+2n^{\frac{-3m^2}{4\|f\|_{\infty}(3+m)}}
\end{equation}
which entails by fixing $m$ and $\varsigma$ large enough
\begin{eqnarray}
T_{6} &\leq& 2(mt_{n})^{4s/(2s+1)}\|f\|^{2}_{W_{2/(1+2s)}} +
t_{n}^{2}\sum_{j\leq J_{n}}\sum_{k}\beta_{jk}^{2} \nonumber \\
&\leq& 2(mt_{n})^{4s/(2s+1)}\|f\|^{2}_{W_{2/(1+2s)}}+
\|f(G^{-1})\|^2_{2}t_{n}^{2}. \nonumber
\end{eqnarray}
Let us look at the term $V_{1}^{"}$
$$V_{1}^{"}=\E\sum_{j\leq
J_{n}}\sum_{k}(w_{jk}^{2}(\hat{\beta}_{jk}-\beta_{jk})^{2}+(1-w_{jk})^{2}\beta_{jk}^{2})I\{|\hat{\beta}_{jk}|>
mt_{n}\}$$
\begin{eqnarray}V_{1}^{"}&=&\E\sum_{j\leq
J_{n}}\sum_{k}(w_{jk}^{2}(\hat{\beta}_{jk}-\beta_{jk})^{2}+(1-w_{jk})^{2}\beta_{jk}^{2})I\{|\hat{\beta}_{jk}|>
mt_{n}\}[I\{ |\beta_{jk}|\leq
mt_{n}/2\}+I\{|\beta_{jk}>mt_{n}/2|\}]\nonumber
\\&=&T_{7}+T_{8}\nonumber\end{eqnarray} for the term $T_{7}$, we use
the Cauchy Scharwz inequality
\begin{eqnarray*}
T_{7}&\leq& \sum_{j\leq J_{n}}\sum_{k}
(\E(\hat{\beta}_{jk}-{\beta}_{jk})^{4})^{1/2}(\p(|\hat{\beta}_{jk}-{\beta}_{jk}|>
mt_{n}/2))^{1/2} \nonumber
\\&+& \sum_{j\leq J_{n}}\sum_{k}
\beta_{jk}^{2}I\{|\hat{\beta}_{jk}|> mt_{n}\}I\{ |\beta_{jk}|\leq
mt_{n}/2\}.\nonumber
\end{eqnarray*}
Furthermore, using inequality $(64)$ p. 1086 in \cite{Kerk2} we get
for $2^j\leq n/\log(n)$
\begin{equation}\label{puissance_4_diff_beta_chap}\E(\hat{\beta}_{jk}-{\beta}_{jk})^{4}\leq
C\frac{1+\|f\|^{4}_{\infty}}{n^2}\end{equation} and by
(\ref{Proba_Beta_jk_chapeau})
\begin{equation*}
\p(|\hat{\beta}_{jk}-{\beta}_{jk}|> mt_{n}/2)\leq
2n^{\frac{-\varsigma^2}{C\|\psi\|^{4}_{4}+\varsigma\|\psi\|^{2}_{\infty}}}+n^{\frac{-m^2}{32(\varsigma+1)}}+2n^{\frac{-3m^2}{16\|f\|_{\infty}(3+m)}}
\end{equation*}
from which follows by fixing again $m$ and $\varsigma$ large enough
\begin{eqnarray} T_{7}&\leq&
C/n.2^{J_{n}}(n^{\frac{-\varsigma^2}{C\|\psi\|^{4}_{4}+\varsigma\|\psi\|^{2}_{\infty}}}+n^{\frac{-m^2}{32(\varsigma+1)}}+2n^{\frac{-3m^2}{16\|f\|_{\infty}(3+m)}}
)^{1/2} +\sum_{j\leq J_{n}}\sum_{k} \beta_{jk}^{2}I\{
|\beta_{jk}|\leq mt_{n}/2\}\nonumber \\&\leq& t_{n}^{2}+
((m/2)t_{n})^{4s/(1+2s)}\|f\|^{2}_{W_{s/(2s+1)}}.
\nonumber\end{eqnarray}
 For the term $T_{8}$
\begin{eqnarray}
T_{8}&=&\E\sum_{j\leq
J_{n}}\sum_{k}(w_{jk}^{2}(\hat{\beta}_{jk}-\beta_{jk})^{2}+(1-w_{jk})^{2}\beta_{jk}^{2})I\{|\hat{\beta}_{jk}|>
mt_{n}\}I\{|\beta_{jk}>mt_{n}/2|\} \nonumber
\\ &\leq&
\frac{4m^{-2/(2s+1)}}{(1-2^{-2/(1+2s)})}\|f\|^{2}_{W_{2/(1+2s)}}({t_{n}})^{4s/(1+2s)}
 \nonumber \\&+&
\E\sum_{j\leq
J_{n}}\sum_{k}(1-w_{jk})^{2}\beta_{jk}^{2}I\{|\hat{\beta}_{jk}|>
mt_{n}\}I\{|\beta_{jk}>mt_{n}/2|\}[I\{|\hat{\beta}_{jk}|\geq|\beta_{jk}/2|\}+I\{|\hat{\beta}_{jk}|<
|\beta_{jk}/2|\}]. \nonumber
\end{eqnarray}
Hereafter we decompose \begin{eqnarray}\E\sum_{j\leq
J_{n}}\sum_{k}(1-w_{jk})^{2}\beta_{jk}^{2})I\{|\hat{\beta}_{jk}|>
mt_{n}\}I\{|\beta_{jk}>mt_{n}/2|\}[I\{|\hat{\beta}_{jk}|\geq|\beta_{jk}/2|\}+I\{|\hat{\beta}_{jk}|<
|\beta_{jk}/2|\}]\nonumber
 \end{eqnarray}
 $$ = T_{8}'+T_{8}^{"}$$
\begin{eqnarray}
T_{8}^{"}\leq \sum_{j\leq J_{n}}\sum_{k}
\beta_{jk}^{2}\p(|\hat{\beta}_{jk}-\beta_{jk}|>mt_{n}/4) \nonumber
\end{eqnarray}
using (\ref{Proba_Beta_jk_chapeau}) we get for $m$ and $\varsigma$
large enough
$$T_{8}^{"}\leq (2n^{\frac{-\varsigma^2}{C\|\psi\|^{4}_{4}+\varsigma\|\psi\|^{2}_{\infty}}}+n^{\frac{-m^2}{128(\varsigma+1)}}+2n^{\frac{-3m^2}{64\|f\|_{\infty}(3+m)}})\sum_{j\leq
J_{n}}\sum_{k}\beta_{jk}^{2}\leq t_{n}^{2}$$\\ as for $T_{8}^{'}$
$$
T_{8}^{'}=\E\sum_{j\leq
J_{n}}\sum_{k}(1-w_{jk})^{2}\beta_{jk}^{2}I\{|\hat{\beta}_{jk}|>
mt_{n}\}I\{|\beta_{jk}|>mt_{n}/2\}I\{|\hat{\beta}_{jk}|\geq|\beta_{jk}|/2\}
$$ using (\ref{lemme2condition3}) we get
\begin{eqnarray*}
T_{8}^{'}&\leq& \E\sum_{j\leq J_{n}}\sum_{k}
K^{2}\beta_{jk}^{2}(\frac{t_{n}}{|\hat{\beta}_{jk}|}+t_{n})^2I\{|\hat{\beta}_{jk}|\geq|\beta_{jk}|/2\}I\{|\beta_{jk}|>mt_{n}/2\}
\nonumber
\\&\leq&
K^{2}\sum_{j\leq
J_{n}}\sum_{k}\beta_{jk}^{2}(\frac{2t_{n}}{|\beta_{jk}|}+t_{n})^{2}I\{|\beta_{jk}|>mt_{n}/2\}
\nonumber
\\ &\leq& 2K^{2}\sum_{j\leq
J_{n}}\sum_{k}\beta_{jk}^{2}(\frac{4t_{n}^{2}}{|\beta_{jk}|^{2}}+t_{n}^2)I\{|\beta_{jk}|>mt_{n}/2\}
\nonumber
\\ &=&
8K^{2}t_{n}^{2}\sum_{j\leq
J_{n}}\sum_{k}I\{|\beta_{jk}|>mt_{n}/2\}+2K^{2}t_{n}^{2}\|f(G^{-1})\|^{2}_{2}
\end{eqnarray*}
using (\ref{weakbesov}) it follows
\begin{eqnarray*}
  T_{8}^{'}&\leq&
8K^{2}t_{n}^{2}{(\frac{mt_{n}}{2})}^{-2/(1+2s)}\frac{2^{2-2/(1+2s)}}{1-2^{-2/(1+2s)}}\|f\|^{2}_{W_{2/(1+2s)}}+2K^{2}t_{n}^{2}\|f(G^{-1})\|^{2}_{2}
\nonumber
\\ &\leq&
32K^2\frac{m^{-2/(1+2s)}}{1-2^{-2/(1+2s)}}t_{n}^{4s/(1+2s)}+2K^{2}t_{n}^{2}\|f(G^{-1})\|^{2}_{2}.\nonumber
\end{eqnarray*}
It remains to bound the bias term $B_{1}$. To this purpose we use
the fact that $f\in B^{s}_{2,\infty}$
$$B_{-1}=\|\sum_{j>J_{n}}\sum_{k=0}^{2^j-1}\beta_{jk}\psi_{j,k}(G(x))\|^{2}_{2} \leq
C2^{-2J_{n}s}=Ct_{n}^{2s} \leq Ct_{n}^{4s/(2s+1)}$$ which completes
the proof.

 \vspace{6mm}
\textbf{Proof of Theorem 2.}\\

In order to prove the Theorem 2., we have to prove that the Bayesian
estimators (\ref{bayesien_beta}) based on Gaussian priors with large
variance (\ref{hyperparametrebigtau}) and (\ref{hyperparametrebigw})
satisfy the conditions of Lemma 2. \\
We will not get into details of the proof because this latter is
identical to the proof of Theorem 3. in \cite{Autin}, with the sole
exception that here, we will place ourselves on the event
$\Omega_{n}^{\delta}$ with $\delta=\varsigma/n$, $\varsigma$ some
positive constant. Indeed, as precised above in section 2.2, a key
observation is that instead of having a deterministic noise
$\eps=1/\sqrt{n}$ like in \cite{Autin}, here we have to deal with a
stochastic noise $\gamma_{jk}^{2}$ which expression is given by
(\ref{bruit_gamma}).

\vspace{6mm}

\end{document}